\documentclass[preprint,11pt]{elsarticle}

\usepackage{amssymb}
\usepackage{amsmath}
\usepackage{amsfonts}
\usepackage{amssymb}
\usepackage{indentfirst,latexsym,bm}
\usepackage{amsthm}
\usepackage{color,xcolor}
\usepackage[all]{xy}
\input{mathrsfs.sty}
\newtheorem{theorem}{Theorem}[section]
\newtheorem{lemma}[theorem]{Lemma}
\newtheorem{corollary}[theorem]{Corollary}

\theoremstyle{definition}
\newtheorem{definition}{Definition}[section]

\theoremstyle{definition}
\newtheorem{example}{Example}[section]

\theoremstyle{remark}
\newtheorem{remark}{Remark}[section]
\theoremstyle{question}

\numberwithin{equation}{section}

\journal{XXX}

\begin{document}

\begin{frontmatter}



\title{A new formula for the weighted Moore-Penrose inverse and its applications}
\author{Qingxiang Xu}
\ead{qingxiang\_xu@126.com}
\address{Department of Mathematics, Shanghai Normal University, Shanghai 200234, PR China}

\begin{abstract}In the general setting of the adjointable operators on Hilbert $C^*$-modules, this paper deals mainly with the weighted Moore-Penrose (briefly weighted M-P) inverse $A^\dag_{MN}$ in the case that the weights $M$ and $N$ are self-adjoint invertible operators, which need not to be positive. A new formula linking $A^\dag_{MN}$ to $A$, $A^\dag$, $M$ and $N$ is derived, in which $A^\dag$ denotes the M-P inverse of $A$.  Based on this formula, some new results on the weighted M-P inverse are obtained. Firstly, it is shown that $A^\dag_{MN}=A^\dag_{ST}$ for some positive definite operators $S$ and $T$. This shows that $A^\dag_{MN}$ is essentially an ordinary weighted M-P inverse. Secondly, some limit formulas for the ordinary weighted M-P inverse originally known for matrices are generalized and improved. Thirdly, it is shown that when $A,M$ and $N$ act on the same Hilbert $C^*$-module, $A^\dag_{MN}$ belongs to the $C^*$-algebra generated by $A$, $M$ and $N$. Finally, some characterizations of the continuity of the weighted M-P inverse are provided.
\end{abstract}


\end{frontmatter}



\section{Introduction}

The weighted Moore-Penrose inverse (briefly weighted M-P inverse) $A^\dag_{MN}$ is closely related to various weighted least squares problems and is dealt with mainly in two cases. The first case is that the weights $M$ and $N$ are positive operators \cite{Corach} (including positive semi-definite matrices \cite{Elden-SINUM,Elden-BIT}), the second case is that $M$ and $N$ are self-adjoint invertible operators \cite{Qin-Xu-Zamani} (including Hermitian nonsingular matrices \cite{Kamaraj-Sivakumar}). The simultaneity of these two cases corresponds to the ordinary weighted M-P inverse, whose weights are positive definite operators (matrices).  Although during the past decades, the weighted M-P inverse and its applications have been intensely studied, there are still some fundamental issues that remain to be unknown.  The purpose of this paper is to put some new insights into these issues in the general setting of the adjointable operators on Hilbert $C^*$-modules.

One of the main issues is investigating formulas for the weighted M-P inverse. Let $A^\dag_{MN}$ be an ordinary weighted M-P inverse. Utilizing the positive definiteness of $M$ and $N$, the following  formula

\begin{equation}\label{not direct relationship}A^\dag_{MN}=N^{-\frac12}(M^\frac12 A N^{-\frac12})^\dag M^{\frac12}\end{equation}
is well-known (see e.g.\,\cite[Lemma~2.3]{SW}) and can be verified directly from the equations stated in \eqref{equ:defn of WPR inverse} below, where $X^\dag$ denotes the M-P inverse of an operator or a matrix $X$.
For matrices $A\in \mathbb{C}^{m\times p}$ and $B\in\mathbb{C}^{n\times p}$ satisfying $\mathcal{N}(A)\cap\mathcal{N}(B)=\{0\}$, the limit formula\footnote{\,The reader should be aware that the notation for the weighted M-P inverse in \cite{Ward} is different from ours.}
\begin{equation}\label{equ:key limit formula}\lim_{t\to 0^+}\left(A^* VA+t B^* W B\right)^{\dag}A^*V=A^\dag_{VU}
\end{equation}
is originally derived in \cite[Theorem~1]{Ward} for  every positive definite matrices $V\in\mathbb{C}^{m\times m}$ and $W\in\mathbb{C}^{n\times n}$, where $U\in\mathbb{C}^{p\times p}$ is a positive definite matrix formulated by $U=A^*VA+B^*WB$. In the special case that $M$ is an identity matrix, some formulas for
$A^\dag_{MN}$ can be found  in \cite[Section~2]{Elden-BIT}, where $N$ is a positive semi-definite matrix. To the best of our knowledge, little has been done on the formulas for $A^\dag_{MN}$ in the case that $M$ and $N$ are self-adjoint invertible operators (matrices). It is significant to find out a formula for $A^\dag_{MN}$ that links $A^\dag_{MN}$ directly to  $A$, $A^\dag$, $M$ and $N$. In this paper, we have managed to do so in the setting of the adjointable operators on Hilbert $C^*$-modules; see \eqref{equ:key formula} in Theorem~\ref{thm:new formula for the wmp} for the details.

Starting from this  newly obtained formula, we are devoted to its applications. Let $H$ and $K$ be Hilbert modules over a $C^*$-algebra $\mathfrak{A}$, let $A$ be an adjointable operator from $H$ to $K$, and let $M$ and $N$ be self-adjoint invertible operators on $K$ and $H$ respectively such that $A^\dag_{MN}$ exists. It is interesting to figure out that
\begin{equation}\label{equ:change weights} A^\dag_{MN}=A^\dag_{S_{A,M}T_{A,N}},\end{equation} where $T_{A,N}\in\mathcal{L}(H)$ and $S_{A,M}\in\mathcal{L}(K)$ are two positive definite operators given by
\begin{align}\label{equ:defn of T}&T_{A,N}=A^\dag A+N(I_H-A^\dag A)N,\\
\label{equ:defn of S}&S_{A,M}=\left[AA^\dag +M^{-1}(I_K-AA^\dag )M^{-1}\right]^{-1}.
\end{align}
This shows that for the above type of weights $M$ and $N$, $A^\dag_{MN}$ is essentially an ordinary weighted M-P inverse.

Another application of formula \eqref{equ:key formula} is the derivation of the limit formulas for the ordinary weighted M-P inverse.
In the matrix case, some limit formulas for the ordinary weighted M-P inverse are derived in  \cite{Ward}. In Section~\ref{sec:limit formula}, we will give a full generalization of these limits formulas. Except that in dealing with Lemma~\ref{lem:decomposition of B to seraration}, our method employed in this section is  quite different from that employed in \cite{Ward}. Furthermore, our main results in this section are derived under the weak condition that $\mathcal{R}(A^*)+\mathcal{R}(B^*)$ is closed, while it is assumed in \cite{Ward} that  $\mathcal{R}(A^*)+\mathcal{R}(B^*)$ is the whole space. For this, two parameters $X$ and $Y$ are introduced in formula \eqref{new form of U with X Y}, so the positive definite weight $U$ defined by this formula is flexible, rather than is fixed as in \cite{Ward}. In addition, as far as we know formula \eqref{wmp separated case} is new even in the matrix case.

Some further applications of formula \eqref{equ:key formula} are concerned with the structure of the $C^*$-algebra generated by $A^\dag_{MN}$ and
the continuity of the weighted M-P inverse. Suppose that $A,M$ and $N$ act on the same Hilbert $C^*$-module such that  $M$ and $N$ are both self-adjoint and invertible. With this restriction on $M$ and $N$, we will show  in Theorem~\ref{thm:where is the C-star alg of WMPI} that whenever  $A^\dag_{MN}$ exists,  $A^\dag_{MN}$  belongs to the $C^*$-algebra generated by $A,M$ and $N$. This enable us to deal with the existence of the weighted M-P inverse by just focusing on the Hilbert space case; see Theorem~\ref{thm:keep wmp invertibility} for the details. The continuity of the M-P inverse of an element in a $C^*$-algebra is dealt with in \cite{Koliha}.
Based on formula \eqref{equ:key formula}, some results obtained in \cite{Koliha} are generalized to the weighted case in the general setting of the adjointable operators on Hilbert $C^*$-modules.

The paper is organized as follows. Section~\ref{sec:1st and 2nd formulas}  is devoted to the derivation of the formulas  \eqref{equ:key formula} and \eqref{equ:change weights}.
The limit formulas for the ordinary weighted M-P inverse are dealt with in Section~\ref{sec:limit formula}. The $C^*$-algebra generated by a weighted M-P inverse and
the continuity of the weighted M-P inverse are concerned with in Sections~\ref{sec:representation of wmp} and \ref{sec:continuity}, respectively.

\section{New formulas for the weighted M-P inverse}\label{sec:1st and 2nd formulas}

Throughout  this paper, $\mathbb{N}$ is the set of positive integers, $\mathbb{C}^{m\times n}$ is the set of $m \times n$ complex matrices, $I_n$ is the identity matrix in $\mathbb{C}^{n\times n}$, $v^T$ denotes the transpose of a vector $v$. Suppose that $\mathfrak{A}$ is a $C^*$-algebra \cite{Murphy,Pedersen},  $H$ and  $K$ are (right) Hilbert module over $\mathfrak{A}$ \cite{Lance,MT,Paschke}.  Let $H\oplus K$ be the induced Hilbert $\mathfrak{A}$-module defined by
$$H\oplus K=\left\{(x,y)^T:x\in H, y\in K\right\},$$ whose $\mathfrak{A}$-valued inner-product is given by
$$\left\langle (x_1,y_1)^T,(x_2,y_2)^T\right\rangle=\langle x_1,x_2\rangle+\langle y_1,y_2\rangle\quad (x_i\in H, y_i\in K, i=1,2).$$ In the  case that $H$ and $K$ are closed submodules of a Hilbert $\mathfrak{A}$-module
such that $H\cap K=\{0\}$, then we use the notation $H\dotplus
K$ to denote the direct sum of $H$ and $K$, that is,
$$H\dotplus
K=\left\{x+y:x\in H, y\in K\right\}.$$
The set of the  adjointable operators from $H$ to $K$ is denoted by $\mathcal{L}(H,K)$ \cite[P.\,8]{Lance}. In case $H=K$, $\mathcal{L}(H,K)$ is abbreviated to $\mathcal{L}(H)$, whose unit (namely, the identity operator on $H$) is denoted  by $I_H$.  Let
$\mathcal{L}(H)_+$ denote the set of positive elements in $\mathcal{L}(H)$.  The notation $T\ge 0$ is also used to indicate that $T$ is an element of $\mathcal{L}(H)_+$.
For each $A\in\mathcal{L}(H,K)$, its adjoint operator, range, null space and restriction on a submodule $H_1$ of $H$ are denoted by $A^*$,
$\mathcal{R}(A)$, $\mathcal{N}(A)$ and $A|_{H_1}$, respectively. It is notable that every element in $\mathcal{L}(H,K)$ is bounded  and  $\mathfrak{A}$-linear, hence $\mathcal{L}(H,K)\subseteq \mathbb{B}(H,K)$, where
$\mathbb{B}(H,K)$ denotes the set of bounded linear operators from $H$ to $K$, with the abbreviation $\mathbb{B}(H)$ for  $\mathbb{B}(H,H)$. It may happen that $\mathcal{L}(H,K)\ne \mathbb{B}(H,K)$ \cite[P.\,8]{Lance}. However, when $\mathfrak{A}=\mathbb{C}$, that is, $H$ and $K$ are Hilbert spaces, then $\mathcal{L}(H,K)=\mathbb{B}(H,K)$. An operator $P\in \mathcal{L}(H)$
is said to be an idempotent if $P=P^2$. If furthermore $P=P^*$, then $P$ is called a  projection. In this case, a unitary operator $U_P: H\to \mathcal{R}(P)\oplus \mathcal{N}(P)$ can be induced as
\begin{equation}\label{equ:the unitray operator induced by a projection P}U_Ph=\big(Ph, (I_H-P)h\big)^T\quad (h\in H)
\end{equation}
such that
\begin{equation}\label{expression of inverse of U P}U_P^{*}\big( (h_1,h_2)^T\big)=h_1+h_2 \quad   \big(h_1\in\mathcal{R}(P), h_2\in\mathcal{N}(P)\big).\end{equation}
By \cite[Sec.~2.2]{Xu-Wei-Gu}  we see that for each $T\in \mathcal{L}(H)$,
 \begin{equation}\label{equ:block matrix T}
U_PTU_P^*=\left(
            \begin{array}{cc}
              PTP|_{\mathcal{R}(P)} & PT(I_H-P)|_{\mathcal{N}(P)} \\
              (I_H-P)TP|_{\mathcal{R}(P)} &(I_H-P)T(I_H-P)|_{\mathcal{N}(P)}\\
            \end{array}
          \right).
\end{equation}
Specifically,
 \begin{equation}\label{equ:block matrix P}
U_PPU_P^*=\left(
            \begin{array}{cc}
              I_{\mathcal{R}(P)} & 0 \\
              0 &0\\
            \end{array}
          \right).
\end{equation}

Unless otherwise specified, throughout the rest of this paper, $\mathfrak{A}$ is a $C^*$-algebra, $H, K$ and $E$  are Hilbert $\mathfrak{A}$-modules.

The main purpose of this section is to derive two new formulas associated with the weighted M-P inverse; see \eqref{equ:key formula} and
\eqref{equ:change weights} for the details. We begin with the definition of the weight on a Hilbert $C^*$-module.
\begin{definition}\label{def:weights}\cite[Sec.~2]{Qin-Xu-Zamani}   An element $M$ of $\mathcal{L}(K)$ is said to be a weight  if $M=M^*$ and $M$ is invertible
in $\mathcal{L}(K)$. If furthermore $M$ is positive, then $M$ is said to be positive definite.
\end{definition}

 Let $M\in \mathcal{L}(K)$ be
a weight. The indefinite inner-product on $K$ induced by $M$ is given by
\begin{equation}\label{quasi-inner-product-induced by M}
\langle x,y\rangle_M=\langle x, My\rangle \quad (x,y\in K).
\end{equation}

\begin{lemma}\label{lem:weighted adjoin operator}{\rm \cite[Lemma~2.2]{Qin-Xu-Zamani}}
Let $M\in\mathcal{L}(K)$ and $N\in\mathcal{L}(H)$ be weights.
Then for every $T\in\mathcal{L}(H,K)$,
\begin{equation*}
\langle Tx,y\rangle_M=\langle x,T^\# y\rangle_N\quad (x\in H, y\in K),
\end{equation*}
where
\begin{equation}\label{equ:weighted adjoint operator-2}
T^\#=N^{-1}T^*M\in\mathcal{L}(K,H).
\end{equation}
\end{lemma}
\begin{remark}\label{rem:def of selfadjiont and projection on HM}The operator $T^\#$ represented as above is called the weighted adjoint operator of $T$. From \eqref{equ:weighted adjoint operator-2}, we see that
$T=T^\#$ iff $(MT)^*=NT$.
\end{remark}

\begin{definition}\label{defn of W-MP inverse}\cite[Sec.~2]{Qin-Xu-Zamani}{\rm
\ Let $M\in \mathcal{L}(K)$ and $N\in
\mathcal{L}(H)$ be weights, and $A\in \mathcal{L}(H,K)$.
The weighted M-P inverse $A^\dag_{MN}$ (if it exists) is the element
$X$ of $\mathcal{L}(K,H)$  satisfying
\begin{equation}\label{equ:defn of WPR inverse}
AXA=A, \quad XAX=X,\quad (MAX)^*=MAX,\quad (NXA)^*=NXA.
\end{equation}
In the special case that $M=I_K$ and $N=I_H$, $A^\dag_{MN}$ is denoted simply by $A^\dag$, which is called the M-P inverse of $A$.
}\end{definition}

Note that if $A^\dag_{MN}$ exists, then it is unique \cite[Lemma~2.3]{Qin-Xu-Zamani}. The existence of the weighted M-P inverse can be clarified as follows.

\begin{lemma}\label{thm:existence of the weighted M-P inverse}{\rm \cite[Theorem~2.4]{Qin-Xu-Zamani}}
Suppose that  $M\in \mathcal{L}(K)$ and $N\in\mathcal{L}(H)$ are weights. Then for every $A\in \mathcal{L}(H,K)$,
$A^\dag_{MN}$ exists if and only if the following conditions are satisfied:
\begin{enumerate}
\item[{\rm (i)}] $\mathcal{R}(A)$ is closed in $K$;
\item[{\rm (ii)}]$\mathcal{R}(AN^{-1}A^*)=\mathcal{R}(A)$;
\item[{\rm (iii)}] $\mathcal{R}(A^*M A)=\mathcal{R}(A^*)$.
\end{enumerate}
\end{lemma}

Specifically, we have the following two lemmas for the M-P inverse.
\begin{lemma}\label{existence of M-P inverse-closedness}{\rm \cite[Theorem~2.2]{Xu-Sheng}}\ For every $A\in \mathcal{L}(H,K)$, $A^\dag$ exists if and
only if $\mathcal{R}(A)$ is closed in $K$.
\end{lemma}

\begin{lemma}\label{lem:orthogonal} {\rm (cf.\,\cite[Theorem 3.2]{Lance} and \cite[Remark 1.1]{Xu-Sheng})}\ For every  $A\in \mathcal{L}(H,K)$,  the closedness of any one of the following sets
implies the closedness of the remaining three sets:
$$\mathcal{R}(A), \quad \mathcal{R}(A^*), \quad \mathcal{R}(AA^*), \quad  \mathcal{R}(A^*A).$$
If $\mathcal{R}(A)$ is closed in $K$, then $\mathcal{R}(A)=\mathcal{R}(AA^*)$,
$\mathcal{R}(A^*)=\mathcal{R}(A^*A)$ and the following  orthogonal
decompositions hold:
\begin{equation*}\label{equ:orthogonal decomposition} H=\mathcal{N}(A)\dotplus
\mathcal{R}(A^*), \quad  K=\mathcal{R}(A)\dotplus \mathcal{N}(A^*).\end{equation*}
\end{lemma}

\begin{remark}If $A\in\mathcal{L}(H,K)$ has the M-P inverse, then $A$ is said to be M-P invertible.
\end{remark}

To derive the main results of this section, we need some additional lemmas obtained in \cite{Qin-Xu-Zamani}.

\begin{lemma}\label{lem: product AA dag and A dag A-one-side fixed}{\rm \cite[Lemma~3.2]{Qin-Xu-Zamani}} Let $M, M_1,M_2\in\mathcal{L}(K)$ and $N, N_1,N_2\in\mathcal{L}(H)$ be weights. If $A,B\in\mathcal{L}(H,K)$ are given such that
$A^\dag_{MN_{1}}$, $A^\dag_{MN_{2}}$, $B^\dag_{M_{1}N}$ and $B^\dag_{M_{2}N}$ exist, then
\begin{equation}\label{equ: product AA dag and A dag A-one-side fixed} AA^\dag_{MN_{1}}=AA^\dag_{MN_{2}},\quad B^\dag_{M_{1}N}B=B^\dag_{M_{2}N}B.
\end{equation}
\end{lemma}

\begin{lemma}\label{lem: A dag M is fixed}{\rm \cite[Theorem~3.4]{Qin-Xu-Zamani}}
Let $M\in\mathcal{L}(K)$ and $N_1,N_2\in\mathcal{L}(H)$ be wei\-ghts.
If $A\in\mathcal{L}(H,K)$ is given such that $A^\dag_{MN_1}$ and $A^\dag_{MN_2}$
exist, then
$A^\dag_{MN_1}=R_{M;N_1,N_2}\cdot A^\dag_{MN_2}$, where
\begin{equation}\label{eqn:defn of R M-N1-N2}
R_{M;N_1,N_2}=A^\dag_{MN_1}A+(I_H-A^\dag_{MN_1}A)N_1^{-1}N_2.
\end{equation}
Furthermore, $R_{M;N_1,N_2}$ is invertible in $\mathcal{L}(H)$ if and only if
\begin{equation*}
\mathcal{R}\Big[(I_H-A^\dag_{MN_1}A)^*N_2(I_H-A^\dag_{MN_1}A)\Big]=\mathcal{R}\big[(I_H-A^\dag_{MN_1}A)^*\big].
\end{equation*}
\end{lemma}

\begin{lemma}\label{lem: A dag N is fixed}{\rm \cite[Theorem~3.5]{Qin-Xu-Zamani}}
Let $M_1,M_2\in\mathcal{L}(K)$ and $N\in\mathcal{L}(H)$ be weights.
If $A\in\mathcal{L}(H,K)$ is given  such that $A^\dag_{M_1N}$ and $A^\dag_{M_2N}$
exist, then
$A^\dag_{M_1N}=A^\dag_{M_2N}\cdot L_{M_1,M_2;N}$, where
\begin{equation}\label{eqn:defn of L M1-M2-N}
L_{M_1,M_2;N}=AA^\dag_{M_1N}+M_2^{-1}M_1(I_K-AA^\dag_{M_1N}).
\end{equation}
Furthermore, $L_{M_1,M_2;N}$ is invertible in $\mathcal{L}(K)$ if and only if
\begin{equation*}
\mathcal{R}\Big[(I_K-AA^\dag_{M_1N})^*\cdot M_1M_2^{-1}M_1\cdot (I_K-AA^\dag_{M_1N})\Big]
=\mathcal{R}\big[(I_K-AA^\dag_{M_1N})^*\big].
\end{equation*}
\end{lemma}

\begin{lemma}\label{lem:existence of different M}{\rm \cite[Theorem~3.10]{Qin-Xu-Zamani}}
Let $M\in\mathcal{L}(K)$ and $N_1,N_2\in\mathcal{L}(H)$ be weights.
If $A\in\mathcal{L}(H,K)$ is given  such that $A^\dag_{MN_1}$ exists and the operator
$R_{M;N_1,N_2}$ defined by \eqref{eqn:defn of R M-N1-N2} is invertible in $\mathcal{L}(H)$.
Then $A^\dag_{MN_2}$ exists such that $A^\dag_{MN_2}=R_{M;N_1,N_2}^{-1}\cdot A^\dag_{MN_1}$.
\end{lemma}

\begin{lemma}\label{lem:existence of different N}{\rm \cite[Theorem~3.11]{Qin-Xu-Zamani}}
Let $M_1,M_2\in\mathcal{L}(K)$ and $N\in\mathcal{L}(H)$ be weights.
If $A\in\mathcal{L}(H,K)$ is given such that $A^\dag_{M_1N}$ exists
and the operator $L_{M_1,M_2;N}$ defined by \eqref{eqn:defn of L M1-M2-N} is invertible in $\mathcal{L}(K)$.
Then $A^\dag_{M_2N}$ exists such that $A^\dag_{M_2N}=A^\dag_{M_1N}\cdot L_{M_1,M_2;N}^{-1}$.
\end{lemma}

\begin{definition}Suppose that $A\in\mathcal{L}(H,K)$ is M-P invertible. For each $X\in\mathcal{L}(H)$ and $Y\in\mathcal{L}(K)$, let
\begin{align}\label{equ:def of R A X}&R_{A,X}=A^\dag A+(I_H-A^\dag A)X,\\
 &\label{equ:def of L A Y} L_{A,Y}=AA^\dag +Y(I_K-AA^\dag).\end{align}
\end{definition}

Now, we are in the position to provide the first main result in this section.

\begin{theorem}\label{thm:new formula for the wmp}Let $M\in\mathcal{L}(K)$ and $N\in\mathcal{L}(H)$ be weights, and $A\in\mathcal{L}(H,K)$. Then the following statements  are  equivalent:
 \begin{enumerate}
 \item[{\rm (i)}] $A^\dag_{MN}$ exists;
 \item[{\rm (ii)}]$A^\dag_{MI_H}$ and  $A^\dag_{I_KN}$ are both existent;
 \item[{\rm (iii)}] $A^\dag$ exists, $R_{A,N}$ and $L_{A,M^{-1}}$ defined by \eqref{equ:def of R A X} and \eqref{equ:def of L A Y} are invertible in $\mathcal{L}(H)$ and $\mathcal{L}(K)$ respectively.
  \end{enumerate}
In each case, $A^\dag_{MN}$ can be represented by
\begin{equation}\label{equ:key formula}A^\dag_{MN}=R_{A,N}^{-1}\cdot  A^\dag \cdot L_{A,M^{-1}}^{-1}.\end{equation}
\end{theorem}
\begin{proof}(i)$\Longleftrightarrow$(ii). If $\mathcal{R}(A)$ is closed in $K$, then by Lemma~\ref{lem:orthogonal}  $\mathcal{R}(A)=\mathcal{R}(AA^*)$ and
$\mathcal{R}(A^*)=\mathcal{R}(A^*A)$. Due to this observation, we see from Lemma~\ref{thm:existence of the weighted M-P inverse} that $A^\dag_{MI_H}$ exists if and only if items (i) and (iii) of Lemma~\ref{thm:existence of the weighted M-P inverse} are both satisfied. Likewise, $A^\dag_{I_KN}$ exists if and only if items (i) and (ii) of Lemma~\ref{thm:existence of the weighted M-P inverse} are both satisfied. This shows the equivalence of (i) and (ii) of this theorem, since $A^\dag_{MN}$ exists if and only if  three items of Lemma~~\ref{thm:existence of the weighted M-P inverse} are all satisfied.

(i)+(ii)$\Longrightarrow$(iii).  Assume that  $A^\dag_{MN}$ exists. By Lemma~\ref{thm:existence of the weighted M-P inverse} $\mathcal{R}(A)$ is closed in $K$, and thus
$A^\dag$ exists according to Lemma~\ref{existence of M-P inverse-closedness}.
For simplicity, we put
$$P_A=AA^\dag,\quad P_{A^*}=A^\dag A.$$
From \eqref{eqn:defn of R M-N1-N2}--\eqref{equ:def of L A Y}, we have
\begin{equation}\label{equ:observation of R and L}R_{A,N}=R_{I_K;I_H,N},\quad L_{A,M^{-1}}=L_{I_K,M;I_H}.\end{equation}
Hence, by Lemmas~\ref{lem: A dag M is fixed} and \ref{lem: A dag N is fixed} we see that
\begin{equation}\label{equ:R and L two-sided}A^\dag=R_{A,N}\cdot A^\dag_{I_KN}=A^\dag_{MI_H}\cdot L_{A,M^{-1}}.
\end{equation}
Furthermore, $R_{A,N}$ is invertible in $\mathcal{L}(H)$ if and only if
\begin{equation}\label{two ranges are the same-01}\mathcal{R}\big[(I_H-P_{A^*})N(I_H-P_{A^*})\big]=\mathcal{R}(I_H-P_{A^*}),\end{equation}
meanwhile $L_{A,M^{-1}}$ is invertible in $\mathcal{L}(K)$ if and only if
\begin{equation}\label{two ranges are the same-02}\mathcal{R}\big[(I_K-P_{A})M^{-1}(I_K-P_{A})\big]=\mathcal{R}(I_H-P_{A}).
\end{equation}

Now, we turn to check the validity of \eqref{two ranges are the same-01} and \eqref{two ranges are the same-02}. By item (ii) of Lemma~\ref{thm:existence of the weighted M-P inverse} and the invertibility of $N$,  we have
\begin{align*}\mathcal{R}(AN^{-1}A^*)=\mathcal{R}(A)=\mathcal{R}(AN^{-1}).
\end{align*}
So for every $x\in H$, there exists $y\in K$ such that $AN^{-1}x=AN^{-1}A^*y$. Put
\begin{align*}&u=(I_H-P_{A^*})N(I_H-P_{A^*})N^{-1}x,\\
&v=-(I_H-P_{A^*})N(I_H-P_{A^*})N^{-1}A^*y.
\end{align*} Utilizing $(I_H-P_{A^*})A^*=0$, we arrive at
$$v=(I_H-P_{A^*})NP_{A^*}N^{-1}A^*y.$$
It follows that
\begin{align*}(I_H-P_{A^*})x=&(I_H-P_{A^*})N\left[(I_H-P_{A^*})+P_{A^*}\right]N^{-1}x\\
=&u+(I_H-P_{A^*})NP_{A^*}N^{-1}A^*y=u+v,
\end{align*}
which clearly leads  to \eqref{two ranges are the same-01}, since $u,v\in \mathcal{R}\big[(I_H-P_{A^*})N(I_H-P_{A^*})\big]$ and $x\in H$ is arbitrary.

Similarly, for every $\xi\in K$, by item (iii) of Lemma~\ref{thm:existence of the weighted M-P inverse} there exists $\eta\in H$ such that
$A^*M\xi=A^*MA\eta$. Let
\begin{align*}&w=(I_K-P_{A})M^{-1}(I_K-P_{A})M\xi,\\
&z=-(I_K-P_{A})M^{-1}(I_K-P_{A})MA\eta.
\end{align*}
Then
\begin{align*}(I_K-P_{A})\xi=&(I_K-P_{A})M^{-1}\left[(I_K-P_{A})+P_{A}\right]M\xi\\
=&w+(I_K-P_{A})M^{-1}(A^\dag)^*A^*M\xi\\
=&w+(I_K-P_{A})M^{-1}(A^\dag)^*A^*MA\eta\\
=&w+(I_K-P_{A})M^{-1}P_{A}MA\eta=w+z,
\end{align*}
which gives \eqref{two ranges are the same-02}, as desired.

Consequently, both $R_{A,N}$ and $L_{A,M^{-1}}$ are invertible, and thus by \eqref{equ: product AA dag and A dag A-one-side fixed} and \eqref{equ:R and L two-sided} we obtain
\begin{align*}A^\dag_{MN}=&A^\dag_{MN}AA^\dag_{MN}=A^\dag_{I_KN}AA^\dag_{MI_H}=R_{A,N}^{-1}A^\dag\cdot A\cdot A^\dag L_{A,M^{-1}}^{-1}\\
=&R_{A,N}^{-1}\cdot A^\dag \cdot L_{A,M^{-1}}^{-1}.
\end{align*}
This shows the validity of \eqref{equ:key formula}.

(iii)$\Longrightarrow$(ii). Suppose that $A^\dag$ exists, and both of $R_{A,N}$ and $L_{A,M^{-1}}$ are invertible. In virtue of \eqref{equ:observation of R and L},  a direct application of Lemmas~\ref{lem:existence of different M} and \ref{lem:existence of different N}
yields the existence of $A^\dag_{I_KN}$ and $A^\dag_{MI_H}$.
\end{proof}

\begin{remark}\label{rem:positive definite weights ensyre}Suppose that $M\in\mathcal{L}(K)$ and $N\in\mathcal{L}(H)$ are both positive definite. If $A\in\mathcal{L}(H,K)$ is given such that $\mathcal{R}(A)$ is closed in $K$,
then $A^\dag$ exists by Lemma~\ref{existence of M-P inverse-closedness}, which  in turn ensures the existence of
$A^\dag_{MN}$ by the positive definiteness of $M$ and $N$ (see e.g.\,\cite[Corollary~2.7]{Qin-Xu-Zamani} or \cite[Theorem~1.3]{Xu-Chen-Song}).
 So, in this case the operators $R_{A,N}$ and $L_{A,M^{-1}}$  defined  by \eqref{equ:def of R A X} and \eqref{equ:def of L A Y}  are both invertible.
\end{remark}

Since every matrix $A$ has a closed range, formula \eqref{equ:key formula} is always valid for positive definite matrices $M$ and $N$.
Based on this formula, some computations of the weighted M-P in the literature can be simplified. We provide such an example as follows.

\begin{example}\cite[Example~5.1]{Xu-Chen-Song} Let $A,M,N\in\mathbb{C}^{4\times 4}$ be given by
$$A=\left(
    \begin{array}{cccc}
      1 & 0 & 1 & -1 \\
      0 & 0 & 1 & 3 \\
      0 & -2 & 0 & 2 \\
      0 & 0 & 0 & 0 \\
    \end{array}\right),\quad M=\left(
                                 \begin{array}{cccc}
                                   2 & 0 & 1 & 0 \\
                                   0 & 1 & 0 & 0 \\
                                   1 & 0 & 1 & 0 \\
                                   0 & 0 & 0 & 1 \\
                                 \end{array}
                               \right),\quad N=\left(
                                                 \begin{array}{cccc}
                                                   2 & 1 & 1 & 0 \\
                                                   1 & 2 & 0 & 0 \\
                                                   1 & 0 & 1 & 0 \\
                                                   0 & 0 & 0 & 1 \\
                                                 \end{array}
                                               \right).
    $$
Then both $M$ and $N$ are positive definite, and
\begin{align*}&A^\dag=\left(
    \begin{array}{cccc}
    11/27      &     1/27     &      2/27    &       0\\
      -4/27    &       7/27   &      -13/27  &         0\\
       4/9     &       2/9    &       -1/18  &         0\\
      -4/27    &       7/27   &        1/54   &        0 \\
        \end{array}
  \right), R_{A,N}=\left(
     \begin{array}{cccc}
      35/27   &       20/27     &     16/27    &       0\\
       2/27   &       32/27     &      4/27    &       0\\
      -2/9    &       -5/9     &       5/9    &        0\\
       2/27   &        5/27    &       4/27    &       1 \\
      \end{array}
   \right),\\
\end{align*}
and $L_{A,M^{-1}}=I_4$.  It follows that
$$A^\dag_{MN}=R_{A,N}^{-1}\cdot A^\dag \cdot L_{A,M^{-1}}^{-1}=\left(
                                    \begin{array}{cccc}
 1/7     &      -2/7    &        3/7     &       0\\
 -3/14    &       5/28  &       -11/28   &        0\\
       9/14     &     13/28   &       -9/28   &        0\\
      -3/14     &      5/28   &        3/28   &        0\\
                                    \end{array}
                                  \right).$$
\end{example}

Our second main result in this section involves the change of wei\-ghts, which indicates that every weighted M-P inverse is actually associated with two positive definite weights.

\begin{theorem}\label{thm:reduce to positive case}Let $M\in\mathcal{L}(K)$ and $N\in\mathcal{L}(H)$ be weights. Suppose that $A\in\mathcal{L}(H,K)$ is given such that $A^\dag_{MN}$ exists. Then \eqref{equ:change weights} is valid, where $T_{A,N}\in\mathcal{L}(H)$ and $S_{A,M}\in\mathcal{L}(K)$ are two positive definite operators given by
\eqref{equ:defn of T} and \eqref{equ:defn of S}, respectively.
\end{theorem}
\begin{proof}For simplicity, we put
$P_A=AA^\dag$ and $P_{A^*}=A^\dag A$ as before. By Theorem~\ref{thm:new formula for the wmp}, the operators $R_{A,N}$ and $L_{A,M^{-1}}$ defined by
\eqref{equ:def of R A X} and \eqref{equ:def of L A Y} are both invertible. Direct computations yield
\begin{equation}\label{equ:change weights++} T_{A,N}=R_{A,N}^*R_{A,N},\quad S_{A,M}=\big(L_{A,M^{-1}}L_{A,M^{-1}}^*\big)^{-1}.\end{equation}
Hence, both $T_{A,N}$ and $S_{A,M}$  are positive definite. We claim that
\begin{equation}\label{equ:N replaced with T}R_{A,N}^{-1}P_{A^*}=R_{A,(T_{A,N})}^{-1}P_{A^*},
\end{equation}
in which $R_{A,(T_{A,N})}$ is defined by \eqref{equ:def of R A X} with $X=T_{A,N}$ therein. That is,
$$R_{A,(T_{A,N})}=P_{A^*}+(I_H-P_{A^*})T_{A,N}.$$
Indeed, \eqref{equ:N replaced with T} is trivially satisfied if $P_{A^*}=0$ or $P_{A^*}=I_H$.  Assume now that $P_{A^*}\ne 0$ and $P_{A^*}\ne I_H$. Let $H_1=\mathcal{R}(P_{A^*}), H_2=\mathcal{N}(P_{A^*})$, and  $U_{P_{A^*}}$ be defined by \eqref{equ:the unitray operator induced by a projection P} such that the projection $P$ therein is replaced with $P_{A^*}$.  Let $I_H$, $I_{H_1}$ and $I_{H_2}$ be denoted simply by $I, I_1$ and $I_2$, respectively. Then
$$U_{P_{A^*}} P_{A^*} U_{P_{A^*}}^*=\left(
                \begin{array}{cc}
                  I_1 & 0 \\
                  0 & 0 \\
                \end{array}
              \right),\quad U_{P_{A^*}}(I-P_{A^*})U_{P_{A^*}}^*=\left(
                                              \begin{array}{cc}
                                                0 & 0 \\
                                                0 & I_2 \\
                                              \end{array}
                                            \right),$$
and
$$U_{P_{A^*}} N U_{P_{A^*}}^*=\left(
                                \begin{array}{cc}
                                  N_{11} & N_{21}^* \\
                                  N_{21} & N_{22} \\
                                \end{array}
                              \right),$$
where
\begin{align*}&N_{11}=P_{A^*}NP_{A^*}|_{H_1},\quad N_{21}=(I-P_{A^*})NP_{A^*}|_{H_1},\\
&N_{22}=(I-P_{A^*})N(I-P_{A^*})|_{H_2}.
\end{align*}
It follows that
\begin{align}U_{P_{A^*}} R_{A,N} U_{P_{A^*}}^*=&U_{P_{A^*}} P_{A^*} U_{P_{A^*}}^*+U_{P_{A^*}}(I-P_{A^*})U_{P_{A^*}}^*\cdot  U_{P_{A^*}} N U_{P_{A^*}}^*\nonumber
\\
=\label{equ:associated MP-invertibility of DP-1}&\left(
                                    \begin{array}{cc}
                                      I_1 & 0 \\
                                      N_{21} & N_{22} \\
                                    \end{array}\right),
                         \end{align}
which implies by the invertibility of $R_{A,N}$ that $N_{22}$ is invertible in $\mathcal{L}(H_2)$. Therefore,
\begin{equation}\label{equ:associated MP-invertibility of DP}U_{P_{A^*}}R_{A,N}^{-1}U_{P_{A^*}}^*=(U_{P_{A^*}} R_{A,N} U_{P_{A^*}}^*)^{-1}=\left(
                                    \begin{array}{cc}
                                      I_1 & 0 \\
                                      -N_{22}^{-1}N_{21} & N_{22}^{-1} \\
                                    \end{array}
                                  \right).
\end{equation}
Consequently,
\begin{align}&\label{equ:useful 2 by 2 expression-01}U_{P_{A^*}}R_{A,N}^{-1}P_{A^*}U_{P_{A^*}}^*=\left(
                                    \begin{array}{cc}
                                      I_1 & 0 \\
                                      -N_{22}^{-1}N_{21} & 0\\
                                    \end{array}
                                  \right),\\
&\label{equ:useful 2 by 2 expression-02}U_{P_{A^*}}R_{A,N}^{-1}(I-P_{A^*})U_{P_{A^*}}^*=\left(
                                    \begin{array}{cc}
                                      0 & 0 \\
                                      0 & N_{22}^{-1}\\
                                    \end{array}
                                  \right).
\end{align}
Meanwhile, direct computations yield
\begin{align*}&R_{A,(T_{A,N})}=P_{A^*}+(I-P_{A^*})N(I-P_{A^*})N,\\
&U_{P_{A^*}}R_{A,(T_{A,N})}U_{P_{A^*}}^*=\left(
                                    \begin{array}{cc}
                                      I_1 & 0 \\
                                      N_{22}N_{21} & N_{22}^2 \\
                                    \end{array}\right),\\
&U_{P_{A^*}}R_{A,(T_{A,N})}^{-1}P_{A^*}U_{P_{A^*}}^*=\left(
                                    \begin{array}{cc}
                                      I_1 & 0 \\
                                      -N_{22}^{-1}N_{21} & 0\\
                                    \end{array}
                                  \right),
\end{align*}
which is combined with \eqref{equ:useful 2 by 2 expression-01} to get \eqref{equ:N replaced with T}. It follows that
\begin{equation}\label{equ:half desired result for positivity}R_{A,N}^{-1}A^\dag=R_{A,N}^{-1}P_{A^*}A^\dag=R_{A,(T_{A,N})}^{-1}P_{A^*}A^\dag=R_{A,(T_{A,N})}^{-1}A^\dag.\end{equation}
Similarly, it can be shown that
\begin{equation}\label{equ:another half desired result for positivity}A^\dag L_{A,M^{-1}}^{-1}=A^\dag L_{A,(S_{A,M})^{-1}}^{-1}.\end{equation}
Actually, it can be derived easily from \eqref{equ:defn of WPR inverse}  that
\begin{equation}\label{dual of WMP inverse}(A^\dag_{MN})^*=(A^*)^\dag_{N^{-1}M^{-1}}.\end{equation} Moreover, according to \eqref{equ:def of R A X} and \eqref{equ:def of L A Y} we have
\begin{align}\label{dual equality-01}&R_{A^*,M^{-1}}=P_A+(I_K-P_A)M^{-1}=(L_{A,M^{-1}})^*,\\
&T_{A^*,M^{-1}}=(R_{A^*,M^{-1}})^*R_{A^*,M^{-1}}=L_{A,M^{-1}}(L_{A,M^{-1}})^*=S_{A,M}^{-1},\nonumber\\
&R_{A^*,(T_{A^*,M^{-1}})}=P_A+(I_K-P_A)S_{A,M}^{-1}=(L_{A,(S_{A,M})^{-1}})^*.\nonumber\end{align}
Replacing $R_{A,N}, A^\dag$ and $R_{A,(T_{A,N})}$ with $R_{A^*,M^{-1}}, (A^*)^\dag$ and $R_{A^*,(T_{A^*,M^{-1}})}$ respectively, it is immediate from \eqref{equ:half desired result for positivity} that
\begin{align*}R_{A^*,M^{-1}}^{-1}(A^*)^\dag=R^{-1}_{A^*,(T_{A^*,M^{-1}})}(A^*)^\dag.
\end{align*}
Hence,
\begin{align*}(L_{A,M^{-1}}^{-1})^*(A^\dag)^*=(L_{A,(S_{A,M})^{-1}}^{-1})^*(A^\dag)^*,
\end{align*}
which gives \eqref{equ:another half desired result for positivity} by taking $*$-operation. Combing \eqref{equ:key formula} with \eqref{equ:half desired result for positivity} and \eqref{equ:another half desired result for positivity} yields
\begin{align*}A^\dag_{MN}=&R_{A,N}^{-1} A^\dag L_{A,M^{-1}}^{-1}=R_{A,N}^{-1} A^\dag\cdot A\cdot A^\dag L_{A,M^{-1}}^{-1}\\
=&R_{A,(T_{A,N})}^{-1} A^\dag\cdot A\cdot A^\dag L_{A,(S_{A,M})^{-1}}^{-1}=A^\dag_{S_{A,M}T_{A,N}}.
\end{align*}
This completes the proof.
\end{proof}

\begin{remark}Suppose that $A^\dag_{MN}$ exists such that $P_{A^*}\ne 0$ and $P_{A^*}\ne I_H$, where $P_{A^*}=A^\dag A$.  It is  useful to determine every positive definite operator $\widetilde{N}\in\mathcal{L}(H)$ that satisfies  \begin{equation}\label{equ:N replaced with T++}R_{A,N}^{-1}P_{A^*}=R_{A,\widetilde{N}}^{-1}P_{A^*},
\end{equation}
since as noted by \eqref{equ:half desired result for positivity} it can be deduced from the above equation that  $R_{A,N}^{-1}A^\dag=R_{A,\widetilde{N}}^{-1}A^\dag$ , and thus
$A^\dag_{I_K N}=A^\dag_{I_K \widetilde{N}}$ according to \eqref{equ:key formula}.
Following the notations in the proof of Theorem~\ref{thm:reduce to positive case}, we have
$$U_{P_{A^*}} \widetilde{N} U_{P_{A^*}}^*=\left(
                                \begin{array}{cc}
                                  \widetilde{N}_{11} & \widetilde{N}_{21}^* \\
                                  \widetilde{N}_{21} & \widetilde{N}_{22} \\
                                \end{array}
                              \right)$$
for some $\widetilde{N}_{ij}\in\mathcal{L}(H_j,H_i)$ ($1\le i,j\le 2$) such that $U_{P_{A^*}} \widetilde{N} U_{P_{A^*}}^*$ is positive definite. Hence, $\widetilde{N}_{22}$ is positive definite and
\begin{align*}U_{P_{A^*}}R_{A,\widetilde{N}}^{-1}P_{A^*}U_{P_{A^*}}^*=\left(
                                    \begin{array}{cc}
                                      I_1 & 0 \\
                                      -\widetilde{N}_{22}^{-1}\widetilde{N}_{21} & 0\\
                                    \end{array}
                                  \right).
\end{align*}
So, \eqref{equ:N replaced with T++} is satisfied if and only if
$$\widetilde{N}_{22}^{-1}\widetilde{N}_{21}=N_{22}^{-1}N_{21}\triangleq N_0.$$
It follows that $\widetilde{N}_{21}=\widetilde{N}_{22}N_0$, which gives
\begin{align}\label{form of widetilde N}U_{P_{A^*}} \widetilde{N} U_{P_{A^*}}^*=&\left(
                                \begin{array}{cc}
                                  \widetilde{N}_{11} & N_0^*\widetilde{N}_{22}\\
                                  \widetilde{N}_{22}N_0 & \widetilde{N}_{22} \\
                                \end{array}
                              \right)\\
=&\left(
    \begin{array}{cc}
      I_1 & N_0^* \\
      0 & I_2 \\
    \end{array}
  \right)
\left(
              \begin{array}{cc}
                \widetilde{N}_{11}-N_0^*\widetilde{N}_{22}N_0 & 0 \\
                0 & \widetilde{N}_{22} \\
              \end{array}
            \right)\left(
                     \begin{array}{cc}
                       I_1 & 0 \\
                       N_0 & I_2 \\
                     \end{array}
                   \right).\nonumber
\end{align}
Therefore, $\widetilde{N}$ can be characterized by \eqref{form of widetilde N} such that both $\widetilde{N}_{22}$ and
$\widetilde{N}_{11}-N_0^*\widetilde{N}_{22}N_0$ are positive definite. Specifically, if we set
$$\widetilde{N}_{22}=N_{22}^2, \quad \widetilde{N}_{21}=N_{22}N_{21},\quad \widetilde{N}_{11}=I_1+N_{21}^*N_{21},$$ then
the operator $T_{A,N}$ defined by  \eqref{equ:defn of T} is derived.

In a similar way, every  positive definite operator $\widetilde{M}\in\mathcal{L}(K)$ satisfying $AA^\dag L_{A,M^{-1}}^{-1}=AA^\dag L_{A,\widetilde{M}^{-1}}^{-1}$  can be characterized. So, there are many positive definite operators $\widetilde{M}\in\mathcal{L}(K)$ and $\widetilde{N}\in\mathcal{L}(H)$ satisfying
$A^\dag_{MN}=A^\dag_{\widetilde{M}\widetilde{N}}$.
\end{remark}

\section{The limit formulas for the ordinary weighted M-P inverse}\label{sec:limit formula}
In this section we restrict our attention to the ordinary weighted M-P inverse, that is, the weights considered here are all positive definite. As noted in Remark~\ref{rem:positive definite weights ensyre},
in this case the M-P invertibility is the same as the weighted  M-P invertibility.
We begin with some known results on the adjointable operators.

\begin{lemma}\label{lem:rang characterization-1}{\rm \cite[Proposition~2.7]{LLX-AIOT}} Let $A\in\mathcal{L}(H,K)$ and $B,C\in\mathcal{L}(E,H)$ be such that $\overline{\mathcal{R}(B)}=\overline{\mathcal{R}(C)}$. Then $\overline{\mathcal{R}(AB)}=\overline{\mathcal{R}(AC)}$.
\end{lemma}

\begin{lemma}\label{lem:Range closure of TT and T} {\rm\cite[Proposition 3.7]{Lance}}
$\overline {\mathcal{R}(T^*T)}=\overline{ \mathcal{R}(T^*)}$ and $\overline {\mathcal{R}(TT^*)}=\overline{ \mathcal{R}(T)}$ for every $T\in\mathcal{L}(H,K)$.
\end{lemma}

\begin{lemma}\label{lem:Range Closure of T alpha and T}{\rm \cite[Lemma~2.3]{XY1}} Suppose that $T\in \mathcal{L}(H)_+$. Then  $\overline{\mathcal{R}(T^{\alpha})}=\overline{\mathcal{R}(T)}$ and $\mathcal{N}(T^{\alpha})=\mathcal{N}(T)$ for every $\alpha>0$.
\end{lemma}

\begin{lemma}\label{lem:closedness is invariant for positive}{\rm \cite[Theorem~2.4]{Vosough-Moslehian-Xu}} For every $T\in\mathcal{L}(H)_+$,  the following statements are  equivalent:
\begin{enumerate}
\item[{\rm (i)}] $\mathcal{R}(T)$ is closed in $H$;
\item[{\rm (ii)}] $\mathcal{R}(T^\alpha)$ is closed in $H$ for all $\alpha>0$;
\item[{\rm (iii)}] $\mathcal{R}(T^\alpha)$ is closed in $H$ for some $\alpha>0$.
\end{enumerate}
\end{lemma}

\begin{lemma}{\rm \cite[Theorem~2.5]{Vosough-Moslehian-Xu}} For every $T\in\mathcal{L}(H)_+$, the following statements are equivalent:
\begin{enumerate}
\item[{\rm (i)}] $\mathcal{R}(T)$ is closed in $H$;
\item[{\rm (ii)}] $\mathcal{R}(T^\alpha)=\mathcal{R}(T)$ for all $\alpha>0$;
\item[{\rm (iii)}] $\mathcal{R}(T^\alpha)=\mathcal{R}(T)$ for some $\alpha\in (0,1)\cup (1,+\infty)$.
\end{enumerate}
\end{lemma}

\begin{lemma}\label{lem:positive operator}{{\rm \cite[Lemma 4.1]{Lance}}} Suppose that $T\in {\mathcal L}(H)$. Then  $T\in {\mathcal L}(H)_+$ if and only if
$\langle Tx, x \rangle\ge 0$ for every $x\in H$.
\end{lemma}

The following lemma will be used in the sequel.

\begin{lemma}\label{lem:C-star elementary}{\rm \cite[Proposition~1.3.5]{Pedersen}} Let $\mathfrak{B}$ be a $C^*$-algebra and $\mathfrak{B}_+$ be the set of positive elements in $\mathfrak{B}$. If $x,y\in \mathfrak{B}_+$ are such that $x\le y$, then $\Vert x\Vert\le \Vert y\Vert$ and $z^*xz\le z^*yz$ for
every $z \in \mathfrak{B}$.
\end{lemma}

Now, we turn to provide some new results.

\begin{lemma}\label{lem:closedness of A+B with parameter} Suppose that $A,B\in\mathcal{L}(H)_+$. Then the following statements are  equivalent:
\begin{enumerate}
\item[{\rm (i)}] $\mathcal{R}(A+B)$ is closed in $H$;
\item[{\rm (ii)}] $\mathcal{R}(A+\lambda B)$ is closed in $H$ for all $\lambda>0$;
\item[{\rm (iii)}] $\mathcal{R}(A+\lambda B)$ is closed in $H$ for some $\lambda>0$.
\end{enumerate}
\end{lemma}
\begin{proof}(i)$\Longrightarrow$(ii). Given $\lambda>0$, let $\alpha=\min\{1,\lambda\}$ and $\beta=\max\{1,\lambda\}$.
It is clear that
$$\alpha(A+B)\le A+\lambda B\le \beta(A+B).$$

Suppose that $y\in\overline{\mathcal{R}\big[(A+\lambda B)^\frac12\big]}$. Then there exists $x_n\in H$ for each $n\in\mathbb{N}$ such that $(A+\lambda B)^\frac12 x_n\to y$.
By Lemmas~\ref{lem:positive operator} and \ref{lem:C-star elementary}, we have
\begin{align*}\left\|(A+B)^\frac12 x_n-(A+B)^\frac12 x_m\right\|^2=&\big\|\langle (A+B)(x_n-x_m), x_n-x_m\rangle\big\|\\
\le&\frac1\alpha \big\|\langle (A+\lambda B)(x_n-x_m), x_n-x_m\rangle\big\|\\
=&\frac1\alpha \left\|(A+\lambda B)^\frac12 x_n-(A+\lambda B)^\frac12 x_m\right\|^2\to 0
\end{align*}
as $n,m\to +\infty$. Hence, there exists $x\in H$ such that $(A+B)^{\frac12} x_n\to (A+B)^{\frac12}x$, since by Lemma~\ref{lem:closedness is invariant for positive} $\mathcal{R}\big[(A+B)^\frac12\big]$ is  closed in $H$. It follows that
\begin{align*}\left\|(A+\lambda B)^\frac12 x_n-(A+\lambda B)^\frac12 x\right\|^2=&\big\|\langle (A+\lambda B)(x_n-x), x_n-x\rangle\big\|\\
\le&\beta\big\|\langle (A+B)(x_n-x), x_n-x\rangle\big\|\\
=&\beta \left\|(A+B)^\frac12 x_n-(A+B)^\frac12 x\right\|^2\to 0.
\end{align*}
Therefore, $y=(A+\lambda B)^\frac12 x$. By the arbitrariness of $y$ in $\overline{\mathcal{R}\big[(A+\lambda B)^\frac12\big]}$, we see that  $\mathcal{R}\big[(A+\lambda B)^\frac12\big]$ is closed in $H$, which leads by Lemma~\ref{lem:closedness is invariant for positive} to the closedness of $\mathcal{R}(A+\lambda B)$.

(ii)$\Longrightarrow$(iii).  It is obvious.

(iii)$\Longrightarrow$(i). Suppose that $\mathcal{R}(A+\lambda B)$ is closed in $H$ for some $\lambda>0$. Then
it can be concluded by the implication (i)$\Longrightarrow$(ii) of this lemma that $\mathcal{R}(A+B)$ is also closed in $H$, since $\lambda B$ is also positive and $A+B$ can be expressed alternatively as $A+\frac{1}{\lambda} (\lambda B)$.
\end{proof}

\begin{lemma}\label{lem:posplam}Let $A,B\in\mathcal{L}(H)_+$. Then the following statements are equivalent:
\begin{enumerate}
\item[{\rm (i)}] $A+B$ is invertible in $\mathcal{L}(H)$;
\item[{\rm (ii)}] $A+\lambda B$ is invertible in $\mathcal{L}(H)$ for all $\lambda>0$;
\item[{\rm (iii)}] $A+\lambda B$ is invertible in $\mathcal{L}(H)$ for some $\lambda>0$.
\end{enumerate}
\end{lemma}
\begin{proof}For every $T\in\mathcal{L}(H)_+$, it is easily seen from Lemma~\ref{lem:orthogonal} that $T$ is invertible in $\mathcal{L}(H)$ if and only if $\mathcal{R}(T)$ is closed in $H$ and $\mathcal{N}(T)=\{0\}$. Since $\mathcal{N}(A+\lambda B)=\mathcal{N}(A)\cap\mathcal{N}(B)=\mathcal{N}(A+B)$, the conclusion is immediate from Lemma~\ref{lem:closedness of A+B with parameter}.
\end{proof}

\begin{lemma}\label{lem:wrt 2 pdos}Let $A\in\mathcal{L}(H,K)$ and $B\in\mathcal{L}(H,E)$ be such that $\mathcal{R}(A^*A+B^*B)$ is closed in $H$. Then for every positive definite operators $V\in\mathcal{L}(K)$ and $W\in\mathcal{L}(E)$, we have
\begin{equation}\label{equ:wrt 2 pdos}\mathcal{R}\left(A^* VA+ B^* W B\right)=\mathcal{R}(A^*A+B^*B)=\mathcal{R}(A^*)+\mathcal{R}(B^*).
\end{equation}
\end{lemma}
\begin{proof}Let $T\in \mathcal{L}(K\oplus E,H)$ be defined by
$T=( A^* V^{\frac12},  B^*W^{\frac12})$, that is,
$$T\left(
     \begin{array}{c}
       x \\
       y \\
     \end{array}
   \right)
=A^* V^{\frac12}x+B^*W^{\frac12}y\quad (x\in K, y\in E).$$ Then $T^*=\begin{pmatrix}V^{\frac12}A\\W^{\frac12}B\end{pmatrix}\in\mathcal{L}(H, K\oplus E)$. Similarly, for $S=(A^*,B^*)\in\mathcal{L}(K\oplus E,H)$ we have $S^*=\begin{pmatrix}A\\B\end{pmatrix}\in\mathcal{L}(H, K\oplus E)$.

It is clear that $\mathcal{R}(T)=\mathcal{R}(S)=\mathcal{R}(A^*)+\mathcal{R}(B^*)$ and $SS^*=A^* A+ B^*  B$. By assumption $\mathcal{R}(SS^*)$ is closed in $H$, so it can be deduced from Lemma~\ref{lem:orthogonal} that $\mathcal{R}(T)=\mathcal{R}(S)=\mathcal{R}(SS^*)$. In particular, $\mathcal{R}(T)$ is closed in $H$. So
we may use Lemma~\ref{lem:orthogonal}  once again to get
$$\mathcal{R}(TT^*)=\mathcal{R}(T)=\mathcal{R}(SS^*),$$
which leads to \eqref{equ:wrt 2 pdos} immediately, since it is obvious that $TT^*=A^* VA+ B^* W B$.
\end{proof}

\begin{lemma}\label{lem:wrt 2 pdosinv}Let $A\in\mathcal{L}(H,K)$ and $B\in\mathcal{L}(H,E)$ be such that $A^*A+B^*B$ is invertible in $\mathcal{L}(H)$. Then for every positive definite operators $V\in\mathcal{L}(K)$ and $W\in\mathcal{L}(E)$,
the operator $A^* VA+ B^* W B$ is invertible in $\mathcal{L}(H)$.
\end{lemma}
\begin{proof}Let $T$ and $S$ be as in the proof of Lemma~\ref{lem:wrt 2 pdos}. By assumption $SS^*$  is invertible in $\mathcal{L}(H)$, so $\mathcal{R}(SS^*)=H$. In particular, $\mathcal{R}(SS^*)$ is closed in $H$. Therefore, by \eqref{equ:wrt 2 pdos}
we have $\mathcal{R}(TT^*)=H$. Evidently,
$$\mathcal{N}(TT^*)=\mathcal{N}(T^*)=\mathcal{N}(A)\cap\mathcal{N}(B)=\mathcal{N}(S^*)=\mathcal{N}(SS^*)=\{0\}.$$
Hence, $TT^*$ is invertible in $\mathcal{L}(H)$.
Since  $TT^*=A^* VA+ B^* W B$, the conclusion follows.
\end{proof}

Our technical lemma in this section reads as follows.
\begin{lemma}\label{lem:limit for the MP inverse-01}Let $A,B\in\mathcal{L}(H)_+$ be such that $A$ is M-P invertible, and $A+B$ is invertible in $\mathcal{L}(H)$. Then the operator  $R_{A,B}$ defined by \eqref{equ:def of R A X}
 is invertible in $\mathcal{L}(H)$ such that
\begin{equation}\label{equ:preparation equality}\lim_{t\to 0^+}\left(A+t B\right)^{-1}A=R_{A,B}^{-1}A^\dag A
\end{equation}
in the operator-norm topology.
\end{lemma}
\begin{proof}Let $P\in\mathcal{L}(H)$ be the projection defined by $P=A^\dag A$. First, we consider the case that $P\ne 0$ and $P\ne I_H$.  Set $H_1=\mathcal{R}(P)$, $H_2=\mathcal{N}(P)$. Let $U_P: H\to H_1\oplus H_2$ be the unitary
operator defined by \eqref{equ:the unitray operator induced by a projection P}, and let $I_H$, $I_{H_1}$ and $I_{H_2}$ be denoted simply by $I, I_1$ and $I_2$, respectively.
Since $A\ge 0$, we have $AA^\dag=A^\dag A=P$, so by \eqref{equ:block matrix T} and \eqref{equ:block matrix P} together with the positivity of $A$ and $B$ it can be deduced that
 \begin{align}\label{equ: unitary representation of A}&U_P PU_P^*=\left(
                            \begin{array}{cc}
                              I_1 & 0 \\
                              0 & 0 \\
                            \end{array}
                          \right),\quad
  U_P A U_P^*=\left(
               \begin{array}{cc}
                 A_1 &0  \\
                0  & 0 \\
               \end{array}
             \right)\in\mathcal{L}(H_1\oplus H_2)_+, \\
 \label{equ: unitary representation of B}&U_P BU_P^*=\left(
               \begin{array}{cc}
                 B_1 & B_2 \\
                 B_2^* & B_3 \\
               \end{array}
             \right)\in\mathcal{L}(H_1\oplus H_2)_+\end{align}
 such that  $A_1=A|_{H_1}\in\mathcal{L}(H_1)_+$ with $A_1^{-1}=A^\dag|_{H_1}$, $B_1=PBP|_{H_1}\in\mathcal{L}(H_1)_+$ and $B_3=(I-P)B(I-P)|_{H_2}\in\mathcal{L}(H_2)_+$.
 By assumption $A+B$ is positive definite, so there exists a positive number $c$ such that $A+B\ge c I$. Hence,
  $$\left(
      \begin{array}{cc}
        A_1+B_1 & B_2 \\
        B_2^* & B_3 \\
      \end{array}
    \right)=U_P(A+B)U_P^*\ge c \left(
                                       \begin{array}{cc}
                                         I_1 &  \\
                                          & I_2 \\
                                       \end{array}
                                     \right),$$
  which implies that $B_3\ge c I_2$. Consequently, $B_3\in\mathcal{L}(H_2)$ is positive definite.

  For every $t>0$, let $T_t=A+t B$, which is positive definite in $\mathcal{L}(H)$  according to Lemma~\ref{lem:posplam}. Hence,
  $U_P T_t U_P^*$ is positive definite in $\mathcal{L}(H_1\oplus H_2)$.  Let $S_t$ be the Schur complement of $U_PT_tU_P^*$ defined by
  $$S_t=tB_3-tB_2^*(A_1+tB_1)^{-1}(tB_2).$$  Direct computation
  yields
  \begin{align*}X_t\cdot  U_PT_tU_P^*\cdot X_t^*=X_t\left(
                               \begin{array}{cc}
                                 A_1+tB_1 & tB_2 \\
                                 tB_2^* & tB_3 \\
                               \end{array}
                             \right)X_t^*
  =\left(
             \begin{array}{cc}
               A_1+tB_1 &  \\
                & S_t \\
             \end{array}
           \right),
  \end{align*}
  where
  $$X_t=\left(
      \begin{array}{cc}
        I_1 & 0 \\
       -tB_2^*(A_1+tB_1)^{-1} & I_2 \\
      \end{array}
    \right).$$
  This shows that $S_t\in\mathcal{L}(H_2)$ is positive definite. It follows that
   \begin{align}
 U_P T_t^{-1}A U_P^*&=X_t^*\left(
                               \begin{array}{cc}
                                 (A_1+tB_1)^{-1} &  \\
                                  & S_t^{-1} \\
                               \end{array}
                             \right)X_t\cdot \left(
                                               \begin{array}{cc}
                                                 A_1 &  \\
                                                  & 0 \\
                                               \end{array}
                                             \right)\nonumber\\
 \label{for convergence-01}&=\left(
               \begin{array}{cc}
                 (A_1+tB_1)^{-1}A_1-t(A_1+tB_1)^{-1}B_2Y(t) & 0 \\
                 Y(t) & 0 \\
               \end{array}
             \right), \end{align}
 where
 \begin{align*}Y(t)=&-tS_t^{-1}B_2^*(A_1+tB_1)^{-1}A_1\\
 =&-\big[B_3-tB_2^*(A_1+tB_1)^{-1}B_2\big]^{-1}B_2^*(A_1+tB_1)^{-1}A_1.
 \end{align*}
 Due to the invertibility of $B_3$ and  $A_1$, we have
 $$\lim_{t\to 0^+}Y(t)=-B_3^{-1}B_2^*$$
 in the operator-norm topology.  This together with  \eqref{for convergence-01} yields
 \begin{equation}\label{equ:matrix form of limit-01}\lim_{t\to 0^+}   U_P T_t^{-1}A U_P^*=\left(
                                              \begin{array}{cc}
                                                I_1 & 0 \\
                                                -B_3^{-1}B_2^* & 0 \\
                                              \end{array}
                                            \right)
                                            \end{equation}
 in the operator-norm topology.

 On the other hand, since $B_3\in\mathcal{L}(H_2)$ is invertible and
 $$U_PR_{A,B}U_P^*=U_P\big[P+(I-P)B\big]U_P^*=\left(
                                   \begin{array}{cc}
                                     I_1 & 0 \\
                                     B_2^* & B_3 \\
                                   \end{array}
                                 \right),
 $$
we see that $R_{A,B}$ is invertible  such that
 $$U_PR_{A,B}^{-1}U_P^*=\left(
                                              \begin{array}{cc}
                                                I_1 & 0 \\
                                                -B_3^{-1}B_2^* & B_3^{-1} \\
                                              \end{array}
                                            \right),$$
 which clearly gives
            \begin{align*}
 U_P R_{A,B}^{-1}PU_P^*=\left(
               \begin{array}{cc}
                 I_1 & 0 \\
                 -B_3^{-1}B_2^*  & 0 \\
               \end{array}
             \right). \end{align*}
 The conclusion is immediate from the above equation together with \eqref{equ:matrix form of limit-01}.

Next, we consider the case that $P=0$ or $P=I$. Observe that \eqref{equ:preparation equality} is trivially satisfied if $P=0$, since in this case we have $A=0$. If $P=I$, then $A$ is invertible in $\mathcal{L}(H)$ such that $A^\dag=A^{-1}$. Hence, both sides of \eqref{equ:preparation equality} are equal to $I$. This completes the proof.
\end{proof}

\begin{remark}To get a deep understanding of the forthcoming formula for the weighted M-P inverse provided in
Theorem~\ref{thm:key limit formula}, it is helpful to give an alternative expression for the right side of \eqref{equ:preparation equality}.
Let $A,B\in\mathcal{L}(H)_+$ be as in Lemma~\ref{lem:limit for the MP inverse-01}. Put $N=A+B$. Then by assumption $N$ is positive definite. Since $A^\dag A=AA^\dag$, we have
$$R_{A,B}=A^\dag A+(I_H-A^\dag A)B=A^\dag A+(I_H-A^\dag A)N=R_{A,N}.$$
It follows from \eqref{equ:preparation equality},  \eqref{equ:key formula} and \eqref{equ: product AA dag and A dag A-one-side fixed} that
\begin{align*}\lim_{t\to 0^+}\left(A+t B\right)^{-1}A=R_{A,N}^{-1}A^\dag A=A^\dag_{I_H N}A=A^\dag_{MN}A
\end{align*}
in the operator-norm topology,
where $M$ is an arbitrary positive definite operator on $H$.
\end{remark}

Based on Theorem~\ref{thm:new formula for the wmp} and Lemma~\ref{lem:limit for the MP inverse-01}, we are now focused on the derivation of some new formulas related to the ordinary weighted M-P inverse.
\begin{corollary}\label{cor:limit for the MP inverse-02}Let $A\in\mathcal{L}(H,K)$ and $B\in\mathcal{L}(H,E)$ be such that $A$ is M-P invertible, and $A^*A+B^*B$ is invertible in $\mathcal{L}(H)$. Then for every positive definite operators $V\in\mathcal{L}(K)$ and $W\in\mathcal{L}(E)$, the operator $R_{A,B^*WB}$ defined by \eqref{equ:def of R A X}
 is invertible in $\mathcal{L}(H)$ such that
\begin{equation*}\label{equ:preparation equality+}\lim_{t\to 0^+}\left(A^*VA+t B^*W B\right)^{-1}A^*VA=R_{A,B^*WB}^{-1}A^\dag A
\end{equation*}
in the operator-norm topology.
\end{corollary}
\begin{proof} By Lemma~\ref{existence of M-P inverse-closedness} $\mathcal{R}(A)$ is closed in $K$, so Lemma~\ref{lem:orthogonal} indicates that
$\mathcal{R}(A^*)$ is closed in $H$ and
$$\mathcal{R}(A^*)=\mathcal{R}(A^* V^{\frac12})=\mathcal{R}\left[(A^* V^{\frac12})(A^* V^{\frac12})^*\right]=\mathcal{R}(A^*VA),$$
which means that $A^*VA$ is M-P invertible such that
\begin{equation}\label{equ: Mp left equal}(A^*VA)^\dag (A^*VA)=A^\dag A.\end{equation} Hence, the desired conclusion can be derived immediately from Lemmas~\ref{lem:wrt 2 pdosinv} and \ref{lem:limit for the MP inverse-01} by replacing $A$ and $B$ with $A^*VA$ and $B^*WB$, respectively.
\end{proof}

\begin{definition}\label{defn:variable U} Suppose that $A\in\mathcal{L}(H,K)$, $B\in\mathcal{L}(H,E)$  and $W\in\mathcal{L}(E)$ are given  such that $W$ is positive definite, and $A^*A+B^*B$ is M-P invertible. Let $\Omega_{A,B,W}$ be the subset of $\mathcal{L}(H)$ consisting of operators $U$ with the form\footnote{\,If $\mathcal{N}(A)\cap\mathcal{N}(B)=\{0\}$, then $U$ is reduced to $A^*XA+B^*WB$.}
\begin{equation}\label{new form of U with X Y}U=A^*XA+B^*WB+YP_{\mathcal{N}(A)\cap\mathcal{N}(B)},\end{equation}
where $P_{\mathcal{N}(A)\cap\mathcal{N}(B)}$ is the projection from $H$ onto $\mathcal{N}(A)\cap\mathcal{N}(B)$, $Y$
is an arbitrary positive definite operator on $\mathcal{N}(A)\cap\mathcal{N}(B)$, and $X\in\mathcal{L}(K)$ is an arbitrary self-adjoint operator  such that $(A^*XA+B^*WB)|_{\mathcal{R}(A^*A+B^*B)}$ is a positive definite operator on $\mathcal{R}(A^*A+B^*B)$.
\end{definition}

\begin{remark} Let $A,B$ and $W$ be given as in Definition~\ref{defn:variable U}.  For every positive definite operator $V\in\mathcal{L}(K)$, it can be deduced from Lemmas~\ref{existence of M-P inverse-closedness} and  \ref{lem:orthogonal} together with \eqref{equ:wrt 2 pdos} that
$$A^*VA+B^*WB+P_{\mathcal{N}(A)\cap\mathcal{N}(B)}\in\Omega_{A,B,W}.$$ Similar reasoning shows that every element $U$ in $\Omega_{A,B,W}$ is positive definite on $H$. So, if furthermore $A$ is M-P invertible, then by Remark~\ref{rem:positive definite weights ensyre}
$R_{A,U}$ is invertible.
\end{remark}

\begin{theorem}\label{thm:key limit formula} Let $A\in\mathcal{L}(H,K)$ and $B\in\mathcal{L}(H,E)$ be such that $A$ and $A^*A+B^*B$ are both M-P invertible. Then for every positive definite operators $V\in\mathcal{L}(K)$ and $W\in\mathcal{L}(E)$, we have
\begin{equation}\label{equ:key limit formula}\lim_{t\to 0^+}\left(A^* VA+t B^* W B\right)^{\dag}A^*V=A^\dag_{VU}
\end{equation}
in the operator-norm topology, where $U$ is an arbitrary element in the set $\Omega_{A,B,W}$ defined by Definition~\ref{defn:variable U}.
\end{theorem}
\begin{proof}If  $A^*A+B^*B=0$, then $A=0$, so \eqref{equ:key limit formula} is trivially satisfied. In what follows we assume that $A^*A+B^*B\ne 0$.

First, we consider the special case that $A^*A+B^*B$ is invertible in $\mathcal{L}(H)$. In this case $\mathcal{N}(A)\cap\mathcal{N}(B)=\{0\}$, which leads by \eqref{equ:def of R A X} and \eqref{new form of U with X Y}  that
$$R_{A,U}=A^\dag A+(I_H-A^\dag A)B^*WB=R_{A,B^*WB}$$
for every $U\in\Omega_{A,B,W}$. In addition, we have
\begin{equation*}A^*VAA^\dag_{VU}=A^*(VAA^\dag_{VU})^*=A^*(A^\dag_{VU})^*A^*V^*=A^*V.
\end{equation*}
It follows from Corollary~\ref{cor:limit for the MP inverse-02}, \eqref{equ:key formula} and \eqref{equ: product AA dag and A dag A-one-side fixed} that
\begin{align*}\lim_{t\to 0^+}\left(A^* VA+t B^* W B\right)^{\dag}A^*V=&\lim_{t\to 0^+}\left(A^* VA+t B^* W B\right)^{\dag}(A^*VA)A^\dag_{VU}\\
=&R_{A,U}^{-1} A^\dag A A^\dag_{VU}=A^\dag_{I_KU} A  A^\dag_{VU}\\=&A^\dag_{VU}A\cdot A^\dag_{VU}=A^\dag_{VU}.
\end{align*}

Next, we consider the general case that $\{0\}\ne H_0\ne H$, where
\begin{equation}\label{defn of H zero}H_0=\mathcal{R}(A^*A+B^*B), \end{equation}
which is closed in $H$ such that $H$ can be decomposed orthogonally as
\begin{equation}\label{equ:decomposition H wrt H0}H=H_0\dotplus H_0^\perp=H_0\dotplus \mathcal{N}(A)\cap \mathcal{N}(B).\end{equation}
Denote by
\begin{equation}\label{Space H 0 induced by A and B}A_0=A^*VA|_{H_0},\quad B_0=B^*WB|_{H_0}.
\end{equation}
From \eqref{equ:wrt 2 pdos} it is clear that
$\mathcal{R}(A_0)\subseteq H_0$ and $\mathcal{R}(B_0)\subseteq H_0$, which mean that $A_0$ and $B_0$ are both in $\mathcal{L}(H_0)_+$.
Let $P_0$ be the projection from $H$ onto $H_0$, and let $U_{P_0}: H\to H_0\oplus H_0^\perp$ be the  unitary defined by \eqref{equ:the unitray operator induced by a projection P} such that  $P$ therein is replaced by $P_0$. Then
\begin{equation}\label{equ:decomposition of A0 B0}U_{P_0} A^*VA U_{P_0}^*=\left(
                                      \begin{array}{cc}
                                        A_0 & 0 \\
                                        0 & 0 \\
                                      \end{array}
                                    \right),\quad U_{P_0}B^*WB U_{P_0}^*=\left(
                                      \begin{array}{cc}
                                        B_0 & 0 \\
                                        0 & 0 \\
                                      \end{array}
                                    \right).
\end{equation}
By the proof of Corollary~\ref{cor:limit for the MP inverse-02}, we see that $A^*VA$ is M-P invertible such that \eqref{equ: Mp left equal} is satisfied.
So, from \eqref{equ:decomposition of A0 B0} it can be concluded that  $A_0$ is M-P invertible in $\mathcal{L}(H_0)$ such that
\begin{align*}U_{P_0}(A^*VA)^\dag A^*VA U_{P_0}^*=\left(
                                      \begin{array}{cc}
                                        A_0^\dag A_0 & 0 \\
                                        0 & 0 \\
                                      \end{array}
                                    \right).
\end{align*}
On the other hand, since $A^\dag A$ is projection and $\mathcal{R}(A^\dag A)=\mathcal{R}(A^*)\subseteq H_0$, we have
\begin{align*}U_{P_0}A^\dag A U_{P_0}^*=
\left(
                \begin{array}{cc}
                  A^\dag A|_{H_0} & 0 \\
                  0 & 0 \\
                \end{array}
              \right).
\end{align*}
This together with \eqref{equ: Mp left equal} yields
\begin{equation}\label{equ:deal with A 0 dag}A_0^\dag A_0= A^\dag A|_{H_0}.
\end{equation}

For each $t>0$, let
\begin{equation}\label{equ:defn of C t}C_t=A^* VA+t B^* W B.\end{equation}
Since $C_t$ can be expressed alternatively as $C_t=A^* VA+B^* (tW) B$, by  \eqref{equ:wrt 2 pdos} we have $\mathcal{R}(C_t)=H_0$, which in turn leads by \eqref{Space H 0 induced by A and B}, \eqref{equ:defn of C t} and  \eqref{equ:decomposition H wrt H0} to   $$\mathcal{R}(A_0+tB_0)=\mathcal{R}\big(C_t|_{H_0}\big)=\mathcal{R}(C_t)=H_0.$$
Evidently, $\mathcal{N}\big(C_t|_{H_0}\big)=H_0\cap \mathcal{N}(A)\cap \mathcal{N}(B)=\{0\}$. So,
$C_t|_{H_0}$ is actually invertible in $\mathcal{L}(H_0)$. It follows from Lemma~\ref{lem:limit for the MP inverse-01} that
\begin{equation}\label{equ:preparation equality++}\lim_{t\to 0^+}\left(C_t|_{H_0}\right)^{-1} A_0=\big[A_0^\dag A_0+(I_{H_0}-A_0^\dag A_0)B_0\big]^{-1}A_0^\dag A_0
\end{equation}
in the operator-norm topology on $\mathcal{L}(H_0)$.

Meanwhile, for each $t>0$ by \eqref{equ:defn of C t} and \eqref{equ:decomposition of A0 B0}  we have
\begin{align*}U_{P_0}C_tU_{P_0}^*=\left(
                                    \begin{array}{cc}
                                      C_t|_{H_0} & 0 \\
                                      0 & 0 \\
                                    \end{array}
                                  \right),
\end{align*}
which implies that
\begin{equation}\label{equ:noname-01}U_{P_0}\cdot C_t^\dag A^*VA \cdot U_{P_0}^*=\left(
                                                              \begin{array}{cc}
                                                               \left(C_t|_{H_0}\right)^{-1} A_0 &  0\\
                                                                0 & 0 \\
                                                              \end{array}
                                                            \right).
\end{equation}

For every $U\in \Omega_{A,B,W}$ given by \eqref{new form of U with X Y}, it is clear that
\begin{align*}R_{A,U}=A^\dag A+(I_H-A^\dag A)B^*WB+(I_H-A^\dag A)Y(I_H-P_0).
\end{align*}
The equation above together with \eqref{equ:deal with A 0 dag} yields
$$U_{P_0} R_{A,U} U_{P_0}^*=\left(
                              \begin{array}{cc}
                                A_0^\dag A_0+(I_{H_0}-A_0^\dag A_0)B_0 &  \\
                                 & Y \\
                              \end{array}
                            \right).
$$
Therefore,
$$U_{P_0} R_{A,U}^{-1} A^\dag A U_{P_0}^*=\left(
                              \begin{array}{cc}
                                \big[A_0^\dag A_0+(I_{H_0}-A_0^\dag A_0)B_0\big]^{-1}A_0^\dag A_0 &  \\
                                 & 0 \\
                              \end{array}
                            \right).$$
We may combine the above equation with \eqref{equ:preparation equality++}, \eqref{equ:noname-01} and  \eqref{equ:key formula}   to obtain
\begin{align*}\lim_{t\to 0^+}C_t^\dag A^*VA=R_{A,U}^{-1} A^\dag A=A^\dag_{I_KU} A=A^\dag_{VU} A.
\end{align*}
Consequently,
\begin{align*}\lim_{t\to 0^+}C_t^\dag A^*V=\lim_{t\to 0^+}C_t^\dag A^*VAA^\dag_{VU}=A^\dag_{VU}.
\end{align*}
This completes the verification of \eqref{equ:key limit formula}.
\end{proof}

\begin{remark}The special case of the preceding theorem can be found in \cite[Theorem~1]{Ward}, where $A,B,V$ and $W$ are matrices such that
$V$ and $W$ are positive definite, and $A^*A+B^*B$ is invertible.  In this case, by Lemma~\ref{lem:wrt 2 pdosinv} $A^*XA+B^*WB$ is positive definite whenever $X\in\mathcal{L}(H)$ is positive definite.  In particular, $A^*VA+B^*WB$ is positive definite.
It is notable that the notation for the weighted M-P inverse in \cite{Ward} is different from ours.  Actually, in our setting the weighted M-P inverse considered in
 \cite[Theorem~1]{Ward} turns out to be $A^\dag_{VU}$, where $U$ is fixed to be $A^*VA+B^*WB$.
 \end{remark}

\begin{remark}Let $A\in\mathcal{L}(H,K)$ and $B\in\mathcal{L}(H,E)$ be such that $A$ is M-P invertible, and $B^*B$ is invertible in $\mathcal{L}(H)$. Then  for every positive definite operators $V\in\mathcal{L}(K)$ and $W\in\mathcal{L}(E)$, one can choose $X=0$ in \eqref{new form of U with X Y}, since
$B^*WB$ is already positive definite in $\mathcal{L}(H)$. This is exactly the case considered in \cite[Corollary~1]{Ward} for matrices.
\end{remark}

We provide a new result as follows.
\begin{theorem}\label{thm:limit for the MP inverse-0001}Let $A,B\in\mathcal{L}(H)_+$ be such that $A$ and $A+B$ are both M-P invertible in $\mathcal{L}(H)$. Then
 $(I_H-A^\dag A)B(I_H-A^\dag A)$ is M-P invertible in $\mathcal{L}(H)$ such that
\begin{equation}\label{equ:lambda tends to infinity-01}\lim_{\lambda\to +\infty}\left(\lambda A+B\right)^\dag B=\big[(I_H-A^\dag A)B(I_H-A^\dag A)\big]^\dag B
\end{equation}
in the operator-norm topology.
\end{theorem}
\begin{proof} First, we consider the special case that $A+B$ is invertible in $\mathcal{L}(H)$. By Lemma~\ref{lem:posplam},    $\lambda A+B$ is invertible
for every $\lambda>0$. If we set $t=\frac{1}{\lambda}$, then
$$(\lambda A+B)^{-1}B=(A+tB)^{-1}(tB)=I_H-(A+tB)^{-1}A.$$
It follows from Lemma~\ref{lem:limit for the MP inverse-01} that
\begin{align*}\label{equ:lambda tends to infinity-02}\lim_{\lambda\to +\infty}\left(\lambda A+B\right)^{-1}B=&I_H-\lim_{t\to 0^+}(A+tB)^{-1}A\\
=&I_H-\big[A^\dag A+(I_H-A^\dag A)B\big]^{-1}A^\dag A\\
=&\big[A^\dag A+(I_H-A^\dag A)B\big]^{-1}(I_H-A^\dag A)B.
\end{align*}
Hence, it is sufficient to prove that
$(I_H-A^\dag A)B(I_H-A^\dag A)$ is M-P invertible in $\mathcal{L}(H)$ such that
\begin{equation}\label{equ:lambda tends to infinity-02}\big[(I_H-A^\dag A)B(I_H-A^\dag A)\big]^\dag=\big[A^\dag A+(I_H-A^\dag A)B\big]^{-1}(I_H-A^\dag A).
\end{equation}
Let $P=A^\dag A$. If $P=0$, then $A+B=B$, so $B$ is invertible, and thus both sides of \eqref{equ:lambda tends to infinity-02} are equal to $B^{-1}$. If $P=I_H$, then
\eqref{equ:lambda tends to infinity-02} is trivially satisfied.

Assume now that $P\ne 0$ and $P\ne I_H$.
Following the notations as in the proof of Lemma~\ref{lem:limit for the MP inverse-01}, we see that $U_P P U_P^*, U_P A U_P^*$ and $U_P B U_P^*$ are given by \eqref{equ: unitary representation of A} and \eqref{equ: unitary representation of B} respectively such that $B_3\in\mathcal{L}(H_2)$ is positive definite. Therefore,
$$U_P(I_H-P)B(I_H-P)U_P^*=\left(
                                          \begin{array}{cc}
                                            0 & 0 \\
                                            0 & B_3 \\
                                          \end{array}
                                        \right),
$$
which is M-P invertible. Consequently, $(I_H-P)B(I_H-P)$ is M-P invertible such that
$$U_P\left[(I_H-P)B(I_H-P)\right]^\dag U_P^*=\left(
                                          \begin{array}{cc}
                                            0 & 0 \\
                                            0 & B_3^{-1} \\
                                          \end{array}
                                        \right).
$$
Meanwhile, if we set $X=\big[A^\dag A+(I_H-A^\dag A)B\big]^{-1}(I_H-A^\dag A)$, then
\begin{align*}U_P XU_P^*=\left(
                           \begin{array}{cc}
                             I_1 & 0 \\
                             B_2^* & B_3 \\
                           \end{array}
                         \right)^{-1}\left(
                                       \begin{array}{cc}
                                         0 & 0 \\
                                         0 & I_2 \\
                                       \end{array}
                                     \right)=\left(
                                          \begin{array}{cc}
                                            0 & 0 \\
                                            0 & B_3^{-1} \\
                                          \end{array}
                                        \right).
\end{align*}
This shows the validity of \eqref{equ:lambda tends to infinity-02}.

Next, we consider the general case that  $A+B$ is M-P invertible in $\mathcal{L}(H)$. Let
 \begin{equation}\label{Space H 0 induced by A and B++}\widetilde{H_0}=\mathcal{R}(A+B), \quad \widetilde{A_0}=A|_{\widetilde{H_0}},\quad \widetilde{B_0}=B|_{\widetilde{H_0}}.\end{equation}
 By assumption $\widetilde{H_0}$ is closed in $H$, so $H$ can be decomposed orthogonally as
 \begin{equation}\label{equ:decomposition H wrt H0++}H=\widetilde{H_0}\dotplus \widetilde{H_0}^\perp=\widetilde{H_0}\dotplus \mathcal{N}(A)\cap \mathcal{N}(B).\end{equation}
 If $\widetilde{H_0}=\{0\}$, then $A=B=0$, so \eqref{equ:lambda tends to infinity-01} is trivially satisfied.  If $\widetilde{H_0}=H$,
then it can be concluded from \eqref{equ:decomposition H wrt H0++} that $A+B$ is invertible in $\mathcal{L}(H)$. Hence, as shown before, \eqref{equ:lambda tends to infinity-01} is also satisfied.

Now, we assume that $\{0\}\ne \widetilde{H_0}\ne H$. In view of Lemma~\ref{lem:wrt 2 pdos}, we have
\begin{align*}\mathcal{R}(A)+\mathcal{R}(B)\subseteq \mathcal{R}(A^\frac12)+\mathcal{R}(B^\frac12)=\mathcal{R}(A+B)\subseteq \mathcal{R}(A)+\mathcal{R}(B).
\end{align*}
This shows that
$$\mathcal{R}(A)+\mathcal{R}(B)=\widetilde{H_0}.$$
Hence, $\mathcal{R}(A)\subseteq \widetilde{H_0}$ and $\mathcal{R}(B)\subseteq \widetilde{H_0}$, which mean that $\widetilde{A_0}$ and $\widetilde{B_0}$ are both in $\mathcal{L}(\widetilde{H_0})_+$ such that $\widetilde{A_0}$ is M-P invertible.
Furthermore, it can be deduced from \eqref{equ:decomposition H wrt H0++} that $\widetilde{A_0}+\widetilde{B_0}$ is invertible in $\mathcal{L}(\widetilde{H_0})$. Therefore,
\begin{equation}\label{equ:defn of C 0}\lim_{\lambda\to +\infty}\left(\lambda \widetilde{A_0}+\widetilde{B_0}\right)^{-1}\widetilde{B_0}=\widetilde{C_0},
\end{equation}
in which $$\widetilde{C_0}=\big[(I_{\widetilde{H_0}}-\widetilde{A_0}^\dag \widetilde{A_0})\widetilde{B_0}(I_{\widetilde{H_0}}-\widetilde{A_0}^\dag \widetilde{A_0})\big]^\dag \widetilde{B_0}.$$

Let $\widetilde{P_0}$ denote the projection from $H$ onto $\widetilde{H_0}$, and let $U_{\widetilde{P_0}}: H\to \widetilde{H_0}\oplus \widetilde{H_0}^\perp$ be the  unitary defined by \eqref{equ:the unitray operator induced by a projection P} such that  $P$ therein is replaced by $\widetilde{P_0}$. Then as shown in the proof of Theorem~\ref{thm:key limit formula}, we have
\begin{equation*}\label{equ:decomposition of A0 B0++}U_{\widetilde{P_0}}A U_{\widetilde{P_0}}^*=\left(
                                      \begin{array}{cc}
                                        \widetilde{A_0} & 0 \\
                                        0 & 0 \\
                                      \end{array}
                                    \right),\ U_{\widetilde{P_0}}A^\dag U_{\widetilde{P_0}}^*=\left(
                                      \begin{array}{cc}
                                        \widetilde{A_0}^\dag  & 0 \\
                                        0 & 0 \\
                                      \end{array}
                                    \right), \ U_{\widetilde{P_0}}B U_{\widetilde{P_0}}^*=\left(
                                      \begin{array}{cc}
                                        \widetilde{B_0} & 0 \\
                                        0 & 0 \\
                                      \end{array}
                                    \right).
\end{equation*}
It follows that
\begin{align*}&U_{\widetilde{P_0}}\left(\lambda A+B\right)^\dag BU_{\widetilde{P_0}}^*=\left(
                                                    \begin{array}{cc}
                                                      \left(\lambda \widetilde{A_0}+\widetilde{B_0}\right)^{-1}\widetilde{B_0} & 0 \\
                                                      0 & 0 \\
                                                    \end{array}
                                                  \right),\\
&U_{\widetilde{P_0}}\big[(I_H-A^\dag A)B(I_H-A^\dag A)\big]^\dag BU_{\widetilde{P_0}}^*=\left(
                                                                  \begin{array}{cc}
                                                                    \widetilde{C_0} & 0 \\
                                                                    0 & 0 \\
                                                                  \end{array}
                                                                \right).
\end{align*}
So, the desired conclusion follows from \eqref{equ:defn of C 0}.
 \end{proof}

\begin{definition}\cite[Definition~1.1 and Remark~1.3]{ELMXZ} Let $A\in\mathcal{L}(H,K)$ and $B\in\mathcal{L}(H,E)$. Then $\big(\mathcal{R}(A^*), \mathcal{R}(B^*)\big)$ is said to be a separated pair if
\begin{equation}\label{equ:isolated condition}\mathcal{R}(A^*)\cap \mathcal{R}(B^*)=\{0\}, \quad \mbox{$\mathcal{R}(A^*)+\mathcal{R}(B^*)$ is closed in $H$}.\end{equation}
\end{definition}

To deal with M-P inverses associated with the separated pairs, we need to make some preparations.
\begin{lemma}\label{lem:overline of R(P+P) and R(P)+R(P) are the same}{\rm \cite[Lemma~2.3 and Remark~5.8]{Luo-Moslehian-Xu}} Let $P,Q\in\mathcal{L}(H)$ be projections. Then
$\overline{\mathcal{R}(P)+\mathcal{R}(Q)}=\overline{\mathcal{R}(P+Q)}.$
Furthermore,  $\mathcal{R}(P)+\mathcal{R}(Q)$ is closed in $H$ if and only if $\mathcal{R}(P+Q)$ is closed in $H$, and in such case $\mathcal{R}(P)+\mathcal{R}(Q)=\mathcal{R}(P+Q)$.
\end{lemma}

\begin{lemma}\label{lem:norm of two projections less than one} Let  $P,Q\in\mathcal{L}(H)$ be projections. Then the following statements are equivalent:
\begin{enumerate}
\item[{\rm (i)}] $\Vert PQ\Vert<1$;
\item[{\rm (ii)}] $\mathcal{R}(P)\cap \mathcal{R}(Q)=\{0\}$ and $\mathcal{R}(P)+\mathcal{R}(Q)$ is closed in $H$;
\item[{\rm (iii)}] $\mathcal{R}(I_H-P)+\mathcal{R}(I_H-Q)=H$;
\item[{\rm (iv)}]$2I_H-P-Q$ is invertible in $\mathcal{L}(H)$.
\end{enumerate}
\end{lemma}
\begin{proof} The equivalence of (i)--(iii) is proved in \cite[Lemma~5.10]{Luo-Moslehian-Xu}.
The implication (iv)$\Longrightarrow$(iii) is obvious.

(iii)+(ii)$\Longrightarrow$(iv). By assumption $\mathcal{R}(I_H-P)+\mathcal{R}(I_H-Q)$ equals the whole space, which is obviously closed. It follows from Lemma~\ref{lem:overline of R(P+P) and R(P)+R(P) are the same} that
$$\mathcal{R}\big[(I_H-P)+(I_H-Q)\big]=H.$$
Hence, $2I_H-P-Q$ is surjective. Furthermore, it is clear that
$$\mathcal{N}(2I_H-P-Q)=\mathcal{R}(P)\cap\mathcal{R}(Q),$$
so by assumption we have $\mathcal{N}(2I_H-P-Q)=\{0\}$.
Therefore, $2I_H-P-Q$ is invertible in $\mathcal{L}(H)$.
\end{proof}

\begin{lemma}\label{lem:closedness technique}Let $A\in\mathcal{L}(H,K)$ and $B\in\mathcal{L}(H,E)$ be such that $\mathcal{R}(A^*)+\mathcal{R}(B^*)$ is closed in $H$. Then the following statements are equivalent:
\begin{enumerate}
\item[{\rm (i)}] $\mathcal{R}(A^*)$ and $\mathcal{R}(B^*)$ are both closed in $H$;
\item[{\rm (ii)}] $\mathcal{R}(A^*A)\cap \mathcal{R}(B^*B)$ is closed in $H$.
\end{enumerate}
\end{lemma}
\begin{proof}Put $Z=A^*A+B^*B$. From the proof of Lemma~\ref{lem:wrt 2 pdos}, we see that $\mathcal{R}(Z)$ is closed in $H$ such that $\mathcal{R}(Z)=\mathcal{R}(A^*)+\mathcal{R}(B^*)$.
Hence, $\mathcal{R}(A^*A)\subseteq \mathcal{R}(A^*)\subseteq \mathcal{R}(Z)$. Similarly, $\mathcal{R}(B^*B)\subseteq \mathcal{R}(Z)$. So, the parallel sum of $A^*A$ and $B^*B$ exists \cite[Theorem~3.4  and Remark~3.5]{Luo-Song-Xu}, which is represented by
$$(A^*A):(B^*B)=A^*A\cdot Z^\dag \cdot B^*B.$$
It follows from \cite[Proposition~4.2]{Luo-Song-Xu} that
$$\mathcal{R}\left[(A^*A):(B^*B)\right]=\mathcal{R}(A^*A)\cap \mathcal{R}(B^*B),$$
and furthermore \cite[Proposition~4.6]{Luo-Song-Xu} indicates that  $\mathcal{R}\left[(A^*A):(B^*B)\right]$ is closed in $H$ if and only if
both $\mathcal{R}(A^*A)$ and $\mathcal{R}(B^*B)$ are closed in $H$. By Lemma~\ref{lem:orthogonal},
the desired conclusion follows immediately.
\end{proof}

\begin{lemma}\label{lem:characterization of separated pair}Suppose that $A\in\mathcal{L}(H,K)$ and $B\in\mathcal{L}(H,E)$ are given such that $\big(\mathcal{R}(A^*), \mathcal{R}(B^*)\big)$ is a separated pair. Then $A$, $B$ and $A^*A+B^*B$ are all M-P invertible such that $2I_H-A^\dag A-B^\dag B$ is invertible in $\mathcal{L}(H)$.
\end{lemma}
\begin{proof} Since $\mathcal{R}(A^*A)\cap \mathcal{R}(B^*B)\subseteq \mathcal{R}(A^*)\cap \mathcal{R}(B^*)=\{0\}$, we have
$\mathcal{R}(A^*A)\cap \mathcal{R}(B^*B)=\{0\}$, so $\mathcal{R}(A^*A)\cap \mathcal{R}(B^*B)$ is closed in $H$.  This, together with the second condition in \eqref{equ:isolated condition}, Lemmas~\ref{lem:closedness technique} and \ref{existence of M-P inverse-closedness}, yields the M-P invertibility of $A$ and $B$.
By assumption $\mathcal{R}(A^*)+\mathcal{R}(B^*)$ is closed in $H$, which gives
$\mathcal{R}(A^*A+B^*B)=\mathcal{R}(A^*)+\mathcal{R}(B^*)$ by using the same technique employed in the proof of Lemma~\ref{lem:wrt 2 pdos}. Hence, by the closedness of
$\mathcal{R}(A^*A+B^*B)$ in $H$, we see that $A^*A+B^*B$ in M-P invertible in $\mathcal{L}(H)$.

Let $P=A^\dag A$ and $Q=B^\dag B$. Then \eqref{equ:isolated condition} can be rephrased as conditions stated in item (ii) of Lemma~\ref{lem:norm of two projections less than one}.
Hence, a direct application of  (ii)$\Longleftrightarrow$(iv) in Lemma~\ref{lem:norm of two projections less than one} yields the invertibility of
$2I_H-A^\dag A-B^\dag B$ in $\mathcal{L}(H)$.
\end{proof}

\begin{remark}\label{rem:an interpretation of separated pair} Let $A\in\mathcal{L}(H,K)$ and $B\in\mathcal{L}(H,E)$.  From Lemma~\ref{lem:norm of two projections less than one} and the proof of Lemma~\ref{lem:characterization of separated pair}, it can be concluded that $\big(\mathcal{R}(A^*), \mathcal{R}(B^*)\big)$ is a separated pair if and only if $A$ and $B$ are both M-P invertible such that
$\|A^\dag A B^\dag B\|<1$.
\end{remark}

We are now in the position to provide a new result as follows.

\begin{theorem}\label{thm:isolated condition}Suppose that $A\in\mathcal{L}(H,K)$ and $B\in\mathcal{L}(H,E)$ are given such that
$\big(\mathcal{R}(A^*), \mathcal{R}(B^*)\big)$ is a separated pair.
Then for every positive definite operators $V\in\mathcal{L}(K)$ and $W\in\mathcal{L}(E)$, we have
\begin{equation}\label{wmp separated case} \left(A^*VA+B^*W B\right)^{\dag}A^*V=(A^*VA)^\dag A^*V-(I_H-A^\dag A)\Pi,\end{equation}
where $\Pi\in\mathcal{L}(K,H)$ is defined by
\begin{equation}\label{defn of C--}\Pi=(2I_H-A^\dag A-B^\dag B)^{-1}(A^*VA)^\dag A^*V.
\end{equation}
\end{theorem}
\begin{proof}By Lemmas~\ref{lem:characterization of separated pair}, \ref{lem:wrt 2 pdos}  and \ref{existence of M-P inverse-closedness} together with  the proof of Corollary~\ref{cor:limit for the MP inverse-02},  we know that
$A$, $B$, $A^*VA+B^*W B$ and $A^*VA$ are all M-P invertible, and $2I_H-A^\dag A-B^\dag B$ is invertible in $\mathcal{L}(H)$. Since $A^*VA$ is positive, we have
$(A^*VA)^\dag A^*VA=A^*VA (A^*VA)^\dag$, which leads by \eqref{equ: Mp left equal} to
\begin{equation}\label{equ:noname001}A^*VA(A^*VA)^\dag A^*=A^\dag A A^*=A^*.\end{equation}
Utilizing \eqref{defn of C--} and $B(I_H-B^\dag B)=0$ yields
\begin{align*}B(A^*VA)^\dag A^*V=B(2I_H-A^\dag A-B^\dag B)\Pi=B(I_H-A^\dag A)\Pi,
\end{align*}
hence $BD=0$, where
\begin{equation}\label{defn of D--} D=(A^*VA)^\dag A^*V-(I_H-A^\dag A)\Pi.\end{equation}
Since $BD=0$ and $A(I_H-A^\dag A)=0$, we may use \eqref{defn of D--} and \eqref{equ:noname001} to get
\begin{align*}\left(A^*VA+B^*W B\right)D=&A^*VAD=A^*VA(A^*VA)^\dag A^*V=A^*V.
\end{align*}
It follows from \eqref{equ:wrt 2 pdos} and \eqref{defn of H zero}  that
\begin{align*}\left(A^*VA+B^*W B\right)^\dag A^*V=&\left(A^*VA+B^*W B\right)^\dag \left(A^*VA+B^*W B\right) D=P_0 D,\end{align*}
where $H_0=\mathcal{R}(A^*)+\mathcal{R}(B^*)$, and $P_0$ is the projection from $H$ onto $H_0$. This together with  \eqref{defn of D--}
indicates that \eqref{wmp separated case} can be simplified as $P_0D=D$.  So, it remains to prove that $\mathcal{R}(D)\subseteq H_0$.

From \eqref{defn of D--}, we see that $D=S-\Pi$, in which $S=(A^*VA)^\dag A^*V+A^\dag A \Pi$ satisfying $\mathcal{R}(S)\subseteq \mathcal{R}(A^*)\subseteq H_0$. Hence, it needs only to show that $\mathcal{R}(\Pi)\subseteq H_0$. Given any $x\in K$, let $u=\Pi x$.
From \eqref{defn of C--}, we have
$$(2I_H-A^\dag A-B^\dag B)u=(A^*VA)^\dag A^*Vx.$$
Consequently,
$$u=\frac12\left[A^\dag A u+B^\dag B u+(A^*VA)^\dag A^*V x\right]\in H_0.$$ So,
 $\mathcal{R}(\Pi)\subseteq H_0$ as desired.
\end{proof}

We may combine Theorems~\ref{thm:key limit formula} and \ref{thm:isolated condition} to get the following corollary.
\begin{corollary}\label{cor:isolated condition}Let $A\in\mathcal{L}(H,K)$, $B\in\mathcal{L}(H,E)$ and $V\in\mathcal{L}(K)$
be given such that $\big(\mathcal{R}(A^*), \mathcal{R}(B^*)\big)$ is a separated pair and $V$ is  positive definite. Then the operator
$\left(A^* VA+t B^* W B\right)^{\dag}A^*V$ equals  $A^\dag_{VU}$, which is independent of the positive number $t$ and the positive definite operator $W\in\mathcal{L}(E)$, where
$U\in\mathcal{L}(H)$ is a positive definite operator given by
\begin{equation}\label{equ:specifical U}U=A^* VA+B^*WB+P_{\mathcal{N}(A)\cap \mathcal{N}(B)}.\end{equation}
\end{corollary}
\begin{proof}Rewrite $A^* VA+t B^* W B$ as $A^* VA+ B^* (tW) B$. The conclusion follows immediately from \eqref{wmp separated case}, \eqref{defn of C--} and \eqref{equ:key limit formula}.
\end{proof}

\begin{remark}The special case of the preceding corollary can be found in \cite[Corollary~2]{Ward} for matrices, where the similar conclusion is derived under the stronger condition that
$$\mathcal{R}(A^*)\cap \mathcal{R}(B^*)=\{0\},\quad \mathcal{R}(A^*)+\mathcal{R}(B^*)=H.$$
\end{remark}

\begin{corollary}Suppose that $A\in\mathcal{L}(H,K)$ and $B\in\mathcal{L}(H,E)$ are  both M-P invertible such that
$$A^\dag A+B^\dag B=I_H.$$
Then for every positive definite operators $V\in\mathcal{L}(K)$ and $W\in\mathcal{L}(E)$, and every positive number $t$, we have
\begin{equation}\label{equ:equiv for weighted MP++} \left(A^* VA+t B^* W B\right)^{-1}A^*V= (A^*VA)^\dag A^*V.\end{equation}
\end{corollary}
\begin{proof}Since $A^\dag A+B^\dag B=I_H$, we have $A^\dag A B^\dag B=0$, so by Remark~\ref{rem:an interpretation of separated pair}  $\big(\mathcal{R}(A^*), \mathcal{R}(B^*)\big)$ is a separated pair. Also, from $A^\dag A+B^\dag B=I_H$
it can be concluded that $A^* VA+B^* (tW) B$ is invertible in $\mathcal{L}(H)$.
Moreover, the operators $\Pi$ and $D$ defined by \eqref{defn of C--} and
\eqref{defn of D--} respectively are both equal to $(A^*VA)^\dag A^*V$. So, the conclusion follows from \eqref{wmp separated case} by replacing $W$ therein with $tW$.
\end{proof}

\begin{remark}The special case of the preceding corollary can be found in \cite[Corollary~2]{Ward} for matrices, where \eqref{equ:equiv for weighted MP++} is derived under the restriction of $V=I_K$ and $W=I_E$.
\end{remark}

\begin{lemma}\label{lem:decomposition of B to seraration} {\rm (cf.\,\cite[Lemma~2]{Ward})} Let $A\in\mathcal{L}(H,K)$ and $B\in\mathcal{L}(H,E)$ be given such that $A$ and $A^*A+B^*B$ are both M-P invertible. Then for every positive definite operators $V\in\mathcal{L}(K)$ and $W\in\mathcal{L}(E)$, there exists a unique decomposition
$$B=B_1+B_2, \quad B_1,B_2\in\mathcal{L}(H,E)$$ satisfying the following  conditions:
 \begin{enumerate}
 \item[{\rm(i)}] $B_2^*WB_1=0$ and  $\mathcal{R}(B_1^*)\subseteq \mathcal{R}(A^*)$;
 \item[{\rm(ii)}] $\big(\mathcal{R}(A^*), \mathcal{R}(B_2^*)\big)$ is a separated pair.
 \end{enumerate}
\end{lemma}
\begin{proof}We provide a proof for the sake of completeness, which is essentially the same as that of \cite[Lemma~2]{Ward}. Let $H_0$ be defined by \eqref{defn of H zero} and  $P_0$ be the projection from $H$ onto $H_0$, and let $U$ be the positive definite operator defined by \eqref{equ:specifical U}. For simplicity, we denote by $Z=A^\dag_{VU}$. Then
$$AZA=A,\quad  (UZA)^*=UZA,$$ and  by \eqref{equ:key limit formula} $\lim\limits_{t\to 0^+}C_t^\dag A^*V=Z$
in the operator-norm topology, where $C_t$ is defined by \eqref{equ:defn of C t}. This shows that
$\mathcal{R}(Z)\subseteq H_0$,
since by \eqref{equ:wrt 2 pdos} $\mathcal{R}(C_t^\dag)=\mathcal{R}(C_t)=H_0$ for every $t>0$.  In view of \eqref{equ:specifical U} and $\mathcal{R}(Z)\subseteq H_0$,  we have $$UZA=A^*VA+B^*WBZA,$$
which gives
\begin{equation}\label{equ:wrt B wrt B}B^*WBZA=A^*Z^*B^*WB,\end{equation} since $UZA$ and $A^*VA$ are both self-adjoint.

Let
\begin{equation}\label{equ:expressions of B1 B2}B_1=B Z A,\quad B_2=B\left(I_H-Z A\right).
\end{equation}
Clearly $B=B_1+B_2$, and  by \eqref{equ:wrt B wrt B} we have
\begin{equation*}\label{equ:wrt B wrt 0}B_2^*WB_1=(I_H-A^*Z^*)B^*WBZA=(I_H-A^*Z^*)A^*\cdot Z^*B^*WB=0,
\end{equation*}
which yields $B_1^*WB_2=0$ by taking $*$-operation.
So
\begin{equation}\label{equ:wrt B  wrt B 2}B^*WB=B_1^*WB_1+B_2^*WB_2.
\end{equation}
Since $B_1^*=(BZA)^*=A^*Z^*B^*$, we have
\begin{equation*}\label{equ:wrt B 1 wrt A}\mathcal{R}\left(B_1^*\right)\subseteq\mathcal{R}\left(A^*\right)\subseteq H_0,
\end{equation*}
which gives
$$\mathcal{R}(B_2^*)=\mathcal{R}(B^*-B_1^*)\subseteq \mathcal{R}(B^*)+\mathcal{R}(B_1^*)\subseteq H_0.$$
Hence
\begin{align*}H_0=\mathcal{R}(A^*)+\mathcal{R}(B^*) \subseteq&\mathcal{R}(A^*)+\mathcal{R}(B_1^*)+\mathcal{R}(B_2^*)\\
=&\mathcal{R}(A^*)+\mathcal{R}(B_2^*) \subseteq H_0,
\end{align*}
which yields $\mathcal{R}(A^*)+\mathcal{R}(B_2^*)= H_0$. So $\mathcal{R}(A^*)+\mathcal{R}(B_2^*)$ is closed in $H$.
If $x\in K$ and $y\in E$ are such that  $A^*x=B_2^*y$, then  $$A^*x=A^*Z^*A^*x=A^*Z^*B_2^*y=A^*\cdot Z^*(I-A^*Z^*)\cdot B^*y=0,$$
hence $\mathcal{R}(A^*)\cap\mathcal{R}(B_2^*)={0}$.
This completes the verification of (i) and (ii) when $B_1$ and $B_2$ are given by \eqref{equ:expressions of B1 B2}.

Conversely, suppose $B$ is decomposed as $B=B_1+B_2$ such that conditions (i) and (ii) are satisfied. Then $(A^*Z^*-I_H)B_1^*=0$, as it is assumed that $\mathcal{R}(B_1^*)\subseteq \mathcal{R}(A^*)$. Hence, $B_1ZA=B_1$. By assumption $B_2^*WB_1=0$, so \eqref{equ:wrt B  wrt B 2} is satisfied. Therefore,
$$B^*WBZA-B_1^*WB_1=(B_1^*WB_1+B_2^*WB_2)ZA-B_1^*WB_1=B_2^*WB_2ZA.$$
Note that $B^*WBZA$ is self-adjoint (see \eqref{equ:wrt B wrt B}), the derivation above shows that $B_2^*WB_2ZA$ turns out to be a subtraction of two self-adjoint operators.
Consequently,  $A^*Z^*B_2^*WB_2=B_2^*WB_2ZA$ and thus $B_2^*WB_2ZA=0$, since by assumption we have  $\mathcal{R}(A^*)\cap\mathcal{R}(B_2^*)={0}$. So
$$\big(W^{\frac{1}{2}}B_2ZA\big)^*W^{\frac{1}{2}}B_2ZA=0.$$
which means that $B_2ZA=0$. Hence,
\begin{align*}&BZA=(B_1+B_2)ZA=B_1ZA=B_1,\\
&B_2=B-B_1=B(I_H-ZA).
\end{align*}
This shows the uniqueness of $B_1$ and $B_2$.
\end{proof}

\begin{theorem}Let $A\in\mathcal{L}(H,K)$ and $B\in\mathcal{L}(H,E)$ be  such that $A$ and $A^*A+B^*B$ are both M-P invertible. Then for every positive definite operators $V\in\mathcal{L}(K)$ and $W,W^\prime\in\mathcal{L}(E)$, we have
\begin{align}\label{wmp change into separated case-01} \lim_{t\to 0^+}\left(A^*VA+t B^*W B\right)^{\dag}A^*V=&\left(A^*VA+B_2^*W^\prime B_2\right)^{\dag}A^*V\\
\label{wmp change into separated case-02}=&(A^*VA)^\dag A^*V-(I_H-A^\dag A)\Pi^\prime\end{align}
in the operator-norm topology,
where $B=B_1+B_2$ is the unique decomposition satisfying conditions (i) and (ii) in Lemma~\ref{lem:decomposition of B to seraration}, and $\Pi^\prime\in\mathcal{L}(K,H)$ is given by
\begin{align}\label{defn of C prime--}&\Pi^\prime=(2I_H-A^\dag A-B_2^\dag B_2)^{-1}(A^*VA)^\dag A^*V.
\end{align}
\end{theorem}
\begin{proof}From the proof of Lemma~\ref{lem:decomposition of B to seraration}, we know that  $B_1=BA^\dag_{VU}A$, where $U$ is a positive definite operator in $\mathcal{L}(H)$  given by \eqref{equ:specifical U}.  Let $C_t$ be defined by \eqref{equ:defn of C t} for each $t>0$, which can be expressed alternatively as
\begin{equation}\label{new decomposition for C t}C_t=A^*VA+tB_1^*WB_1+tB_2^*WB_2=A^*V_tA+B_2^*(tW)B_2,\end{equation}
in which
\begin{equation*}\label{equ:defn of V t} V_t=V+t\left(BA^\dag_{VU}\right)^*W(BA^\dag_{VU}).\end{equation*}
Note that $\big(\mathcal{R}(A^*),\mathcal{R}(B_2^*)\big)$ is a separated pair, and the operator $V_t$ given as above is positive definite. So, by \eqref{new decomposition for C t} and Theorem~\ref{thm:isolated condition} we have
\begin{equation}\label{wmp separated case with parameter} C_t^\dag A^*V_t=D_t,\end{equation}
where $\Pi_t,D_t\in\mathcal{L}(K,H)$ are given by
\begin{align}\label{defn of pi t--}&\Pi_t=(2I_H-A^\dag A-B_2^\dag B_2)^{-1}(A^*V_tA)^\dag A^*V_t,\\
\label{defn of D t--}&D_t=(A^*V_tA)^\dag A^*V_t-(I_H-A^\dag A)\Pi_t.\end{align}

Now, let $P$, $H_1$, $H_2$ and $U_P$ be defined as in the proof of Lemma~\ref{lem:limit for the MP inverse-01}. Suppose firstly that $P\ne 0$ and $P\ne I_H$.
In this case, it is clear that
\begin{align*}&U_P A^*V_tA U_P^*=\left(
                     \begin{array}{cc}
                       A^*V_tA|_{H_1} & 0 \\
                       0 & 0 \\
                     \end{array}
                   \right),\\
&U_P (A^*V_tA)^\dag  U_P^*=\left(
                     \begin{array}{cc}
                       \left(A^*V_tA|_{H_1}\right)^{-1} & 0 \\
                      0 & 0 \\
                     \end{array}
                   \right).
\end{align*}
So, due to $\lim\limits_{t\to 0^+}\|V_t-V\|=0$ we have
$$\lim_{t\to 0^+}\big\| A^*V_tA|_{H_1}- A^*VA|_{H_1}\big\|=\lim_{t\to 0^+}\big\| A^*V_tA- A^*VA\big\|=0,$$
which means that
$$\lim_{t\to 0^+}\big\|( A^*V_tA)^\dag- (A^*VA)^\dag\big\|=\lim_{t\to 0^+}\left\| \left(A^*V_tA|_{H_1}\big)^{-1}-\big( A^*VA|_{H_1}\right)^{-1}\right\|=0.$$
It follows from \eqref{defn of pi t--}, \eqref{defn of C prime--} and \eqref{defn of D t--}  that
$$\lim_{t\to 0^+}\|\Pi_t-\Pi^\prime\|=0,\quad  \lim_{t\to 0^\dag}\|D_t-D^\prime\|=0,$$
where  $\Pi^\prime$ is defined by \eqref{defn of C prime--} and
\begin{equation*}D^\prime=(A^*VA)^\dag A^*V-(I_H-A^\dag A)\Pi^\prime,\end{equation*}
which is equal to $\left(A^*VA+B_2^*W^\prime B_2\right)^{\dag}A^*V$ whenever $W^\prime\in\mathcal{L}(E)$ is an arbitrary positive definite operator, as pointed out by \eqref{wmp separated case}.
Hence, by \eqref{wmp separated case with parameter} we obtain
$$\lim_{t\to 0^\dag} C_t^\dag A^*V=\lim_{t\to 0^\dag} C_t^\dag A^*V_t\cdot (V_t^{-1}V)=D^\prime \cdot I_K=D^\prime.$$
This shows the validity of \eqref{wmp change into separated case-01} and \eqref{wmp change into separated case-02}.

Finally, we consider the case that $P=0$ or $P=I_H$. If $P=0$, then $A=0$, so \eqref{wmp change into separated case-01} and \eqref{wmp change into separated case-02} are trivially satisfied. If $P=I_H$, then $H_1=H$, $A^*V_tA$ $(t>0)$ and $A^*VA$ are invertible. Similar reasoning shows that  \eqref{wmp change into separated case-01} and \eqref{wmp change into separated case-02} are also satisfied.
\end{proof}

\begin{remark}The limit formula \eqref{wmp change into separated case-01} is derived in \cite[Theorem~2]{Ward} for matrices in the special case
that $\mathcal{R}(A^*)+\mathcal{R}(B^*)=H$, while formula \eqref{wmp change into separated case-02} is newly obtained.
\end{remark}

\section{The $C^*$-algebra generated by a weighted M-P inverse}\label{sec:representation of wmp}

In this section, we restrict our attention to the case that $H=K$. Suppose that $M,N\in\mathcal{L}(H)$ are two weights, which are not necessary to be positive definite. Based on formula \eqref{equ:key formula}, we study mainly  the $C^*$-algebra generated by a weighted M-P inverse $A^\dag_{MN}$.

\begin{definition}For every $C^*$-algebra $\mathfrak{B}$ and its subset $\mathfrak{F}$, let $C^*\mathfrak{F}$ be the $C^*$-subalgebra of $\mathfrak{B}$ generated by $\mathfrak{F}$.
\end{definition}
For each $T\in\mathcal{L}(H)$, by definition  $C^*\{T,I_H\}$ denotes the $C^*$-subalgebra of $\mathcal{L}(H)$ generated by operators $T$ and $I_H$.

\begin{lemma}\label{lem:inverse is in}{\rm \cite[Theorem~2.1.11]{Murphy}} Let $\mathfrak{C}$ be a $C^*$-subalgebra of a unital $C^*$-algebra $\mathfrak{B}$ containing the unit of $\mathfrak{B}$. Then for every $c\in \mathfrak{C}$, $c$ is invertible in $\mathfrak{B}$ if and only if it is invertible in $\mathfrak{C}$.
\end{lemma}

\begin{lemma}\label{lem: A dag contained sub alg} Suppose that $T\in\mathcal{L}(H)$ is M-P invertible. Then $T^\dag$ belongs to the $C^*$-subalgebra of $\mathcal{L}(H)$ generated by $T$.
\end{lemma}
\begin{proof}Let $K=E=H$, $A=T$, $B=V=W=I_H$, and let $X=0$ and $Y=0$. Then by \eqref{new form of U with X Y} and \eqref{equ:key limit formula}, we have
$U=I_H\in\Omega_{A,B,W}$ and
\begin{equation}\label{equ:limit formula for MP}\lim_{t\to 0^+} (T^*T+t I_H)^{-1}T^*=T^\dag
\end{equation}
in the operator-norm topology. Hence, by \eqref{equ:limit formula for MP} and Lemma~\ref{lem:inverse is in} we see that
$$T^\dag\in C^*\{T,I_H\}\cdot C^*\{T\}=C^*\{T\}.\qedhere$$
\end{proof}

\begin{example}Given an idempotent $Q\in\mathcal{L}(H)$, let $m(Q)$ be its matched projection introduced in \cite[Sec.~2]{TXF02} as
\begin{equation*}m(Q)=\frac12\big(|Q^*|+Q^*\big)|Q^*|^\dag\big(|Q^*|+I_H\big)^{-1}\big(|Q^*|+Q\big).
\end{equation*}
According to Lemmas~\ref{lem:inverse is in} and \ref{lem: A dag contained sub alg}, we have
$$m(Q)\in C^*\{Q\} \cdot C^*\{Q, I_H\}\cdot C^*\{Q\}=C^*\{Q\}.\qedhere$$
\end{example}

\begin{theorem}\label{thm:where is the C-star alg of WMPI}Suppose that  $M,N\in\mathcal{L}(H)$ are two weights in the sense of Definition~\ref{def:weights}. If $A\in\mathcal{L}(H)$ is given such that $A^\dag_{MN}$ exists. Then $A^\dag_{MN}$ belongs to the $C^*$-subalgebra of $\mathcal{L}(H)$  generated by $A,M$ and $N$.
\end{theorem}
\begin{proof}By Theorem~\ref{thm:new formula for the wmp} $A^\dag$ is M-P invertible, and the operators $R_{A,N}$ and $L_{A,M^{-1}}$ defined by \eqref{equ:def of R A X} and \eqref{equ:def of L A Y} are both invertible such that \eqref{equ:key formula} is satisfied. From Lemma~\ref{lem: A dag contained sub alg} we have $A^\dag\in C^*\{A\}$, which is combined with  Lemma~\ref{lem:inverse is in} to conclude that
$$R_{A,N}^{-1}\in C^*\{A, N, I_H\}, \quad  L_{A,M^{-1}}^{-1}\in C^*\{A, M, I_H\}.$$
Hence, it can be derived from \eqref{equ:key formula} that
$$A^\dag_{MN}\in C^*\{A, N, I_H\} \cdot C^*\{A\}\cdot C^*\{A, M, I_H\}\subseteq C^*\{A,M,N\}.\qedhere$$
 \end{proof}

Recall that a pair $(\pi, X)$ is said to be a representation of a $C^*$-algebra $\mathfrak{B}$, if $X$ is a Hilbert space and $\pi: \mathfrak{B}\to \mathbb{B}(X)$ is a $C^*$-morphism.

\begin{theorem}\label{thm:keep wmp invertibility}
 Let $(\pi, X)$ be a faithful unital representation  of $\mathcal{L}(H)$. Then for every weights $M,N\in\mathcal{L}(H)$ and every operator $A\in \mathcal{L}(H)$,
 $A^\dag_{MN}$ exists  if and only if $\pi(A)^\dag_{\pi(M)\pi(N)}$ exists. In such case,
$\pi(A)^\dag_{\pi(M)\pi(N)}=\pi(A^\dag_{MN})$.
 \end{theorem}
\begin{proof} Since $(\pi, X)$ is a unital representation, both $\pi(M)$ and $\pi(N)$ are weights in $\mathbb{B}(X)$. To simplify the notation, we put $\mathfrak{B}=C^*\{A,M,N\}$.  It is well-known that  $\|\pi\|\le 1$ and $\pi(\mathfrak{B})$ is a $C^*$-subalgebra of $\mathbb{B}(X)$ (see e.g.\,\cite[Theorem~1.5.7]{Pedersen}), which lead clearly to
$\pi(\mathfrak{B})\subseteq C^*\{\pi(A),\pi(M),\pi(N)\}$ and $\pi(\mathfrak{B})\supseteq C^*\{\pi(A),\pi(M),\pi(N)\}$, respectively. Therefore,
\begin{equation}\label{equ:two c star algs are equal-weighted}\pi(\mathfrak{B})=C^*\{\pi(A),\pi(M),\pi(N)\}.\end{equation}

Suppose now that $A^\dag_{MN}$ exists. Let  $A^\dag_{MN}$ be denoted simply by $Z$. Then from \eqref{equ:defn of WPR inverse}, we have
\begin{align}\label{equ: Z as wmp}AZA=A, \quad ZAZ=Z,\quad \big(MAZ\big)^*=MAZ, \quad \big(NZA\big)^*=NZA,
\end{align}
which gives
\begin{align}\label{wmp 12 equations for p--i} &\pi(A)\pi(Z)\pi(A)=\pi(A), \quad \pi(Z)\pi(A)\pi(Z)=\pi(Z),\\
\label{wmp 3 equation for pi--} &\big[\pi(M)\pi(A)\pi(Z)\big]^*=\pi(M)\pi(A)\pi(Z), \\
\label{wmp 4 equation for pi--}&\big[\pi(N)\pi(Z)\pi(A)\big]^*=\pi(N)\pi(Z)\pi(A).
\end{align}
  This shows that $\pi(A)^\dag_{\pi(M)\pi(N)}$  exists such that $\pi(A)^\dag_{\pi(M)\pi(N)}=\pi\big(Z)$.

Conversely, assume that $\pi(A)^\dag_{\pi(M)\pi(N)}$ exists. By Theorem~\ref{thm:where is the C-star alg of WMPI},
we have $\pi(A)^\dag_{\pi(M)\pi(N)}\in C^*\{\pi(A),\pi(M),\pi(N)\}$. In view of \eqref{equ:two c star algs are equal-weighted}, there exists $Z\in\mathfrak{B}$ such that
$\pi(A)^\dag_{\pi(M)\pi(N)}=\pi(Z)$. So \eqref{wmp 12 equations for p--i}--\eqref{wmp 4 equation for pi--} are all satisfied. As a result,
\eqref{equ: Z as wmp} can be derived from the faithfulness of $\pi$. This shows that $Z=A^\dag_{MN}$.
\end{proof}

\section{The continuity of the weighted M-P inverse}\label{sec:continuity}

It is notable that Section~\ref{sec:representation of wmp} deals merely with $A^\dag_{MN}$ in the special case that $A,M$ and $N$ act on the same Hilbert $C^*$-module. To
handle norm estimations for the weighted M-P inverse
in the general case, we need the following auxiliary lemma.

\begin{lemma}\label{lem:self-adjointable induced by rho}Given an operator $A\in\mathcal{L}(H,K)$, and  weights $M\in\mathcal{L}(K)$ and $N\in\mathcal{L}(H)$, let $\rho(A),T\in\mathcal{L}(K\oplus H)$ be defined by
\begin{equation*}\label{equ:defn of the morphism rho}\rho(A)=\left(
                          \begin{array}{cc}
                            0 & A \\
                            A^* & 0 \\
                          \end{array}
                        \right),\quad T=\left(
                                          \begin{array}{cc}
                                            M & 0 \\
                                            0 & N^{-1} \\
                                          \end{array}
                                        \right).
\end{equation*}
Then $A^\dag_{MN}$    exists if and only if  $\rho(A)^\dag_{TT^{-1}}$ exists. In such case, we have
\begin{equation*}\label{equ:basic property of rho A}\rho(A)^\dag_{TT^{-1}}=\left(
                               \begin{array}{cc}
                                 0 & (A^\dag_{MN})^* \\
                                 A^\dag_{MN} & 0\\
                               \end{array}
                             \right).
\end{equation*}
\end{lemma}
\begin{proof} A direct use of \eqref{equ:defn of WPR inverse} yields the desired conclusion.
\end{proof}

Now, we are in the position to provide our first result on the continuity of the weighted M-P inverse.

\begin{theorem}Let operators $A,A_n\in\mathcal{L}(H,K)\setminus\{0\}$, and weights $M_n,M\in\mathcal{L}(K)$, $N_n, N\in\mathcal{L}(H)$ be given such that
 $A_n\to A, M_n\to M$ and $N_n\to N$ in the operator-norm topology. Suppose that $A^\dag_{MN}$ exists, and $(A_n)^\dag_{M_nN_n}$ exists for every $n\in \mathbb{N}$.
 Then the following conditions are equivalent:
 \begin{enumerate}
\item[{\rm (i)}] $\lim\limits_{n\to\infty}\big\|(A_n)^\dag_{M_nN_n}-A^\dag_{MN}\big\|=0$;
\item[{\rm (ii)}]$\sup\limits_n\big\{\big\|(A_n)^\dag_{M_nN_n}\big\|\big\}<+\infty$;
\item[{\rm (iii)}]$\sup\limits_n\big\{\|A_n^\dag\|\big\}<+\infty$;
\item[{\rm (iv)}] $\lim\limits_{n\to \infty}\big\|A_n^\dag A_n-A^\dag A\big\|=0$;
\item[{\rm (v)}] $\lim\limits_{n\to \infty}\big\| A_n A_n^\dag -A A^\dag \big\|=0$;
\item[{\rm (vi)}] $\lim\limits_{n\to \infty}\big\|A_n^\dag-A^\dag\big\|=0$.
\end{enumerate}
 \end{theorem}
 \begin{proof}In view of Lemma~\ref{lem:self-adjointable induced by rho} and the observation of
 $\|\rho(T)\|=\|T\|$ for every $T\in\mathcal{L}(H,K)$ \cite[Lemma~4.1]{Xu-Wei-Gu}, we may as well assume that $H=K$. In this case, all the operators considered are elements in the $C^*$-algebra $\mathcal{L}(H)$, so the equivalence of (iii)--(vi) is established in \cite[Theorem~1.6]{Koliha}. The implication (i)$\Longrightarrow$(ii) is clear. Hence, it remains only to prove the implications (vi)$\Longrightarrow$(i) and (ii)$\Longrightarrow$(iii).

 (vi)$\Longrightarrow$(i). It follows from Theorem~\ref{thm:new formula for the wmp} that
 \begin{equation}\label{equ:relationships repeated for n}(A_n)^\dag_{M_nN_n}=R_{A_n,N_n}^{-1}\cdot A_n^\dag\cdot L_{A_n,M_n^{-1}}^{-1},\quad A^\dag_{MN}=R_{A,N}^{-1}\cdot A^\dag\cdot L_{A,M^{-1}}^{-1},\end{equation}
 where $R_{A_n,N_n}, R_{A,N}$ and $L_{A_n,M_n^{-1}}, L_{A,M^{-1}}$ are defined by \eqref{equ:def of R A X} and \eqref{equ:def of L A Y}, respectively.
By the assumption, it is clear that
$R_{A_n, N_n}\to R_{A,N}$ and $L_{A_n, M_n}\to L_{A,M^{-1}}$, which in turn lead to
$R_{A_n, N_n}^{-1}\to R_{A,N}^{-1}$ and $L_{A_n, M_n^{-1}}^{-1}\to L_{A,M^{-1}}^{-1}$. Utilizing \eqref{equ:relationships repeated for n} yields $(A_n)^\dag_{M_nN_n}\to A^\dag_{MN}$.

 (ii)$\Longrightarrow$(iii). From the assumption together with the observation of $M_n^{-1}\to M^{-1}$ and $N_n\to N$ in the operator-norm topology, we have
  $$\alpha:=\sup_n\left\{\left\|M_n^{-1}\right\|,\|N_n\|,\left\|(A_n)^\dag_{M_nN_n}\right\|\right\}<+\infty.$$
So for every $n\in \mathbb{N}$,
\begin{align*}&\|R_{A_n,N_n}\|\le \|A_n^\dag A_n\|+\|(I_H-A_n^\dag A_n)N_n\|\le 1+\|N_n\|\le 1+\alpha,\\
&\|L_{A_n,M_n^{-1}}\|\le \|A_nA_n^\dag \|+\left\|M_n^{-1}(I_H-A_n A_n^\dag)\right\|\le 1+\left\|M_n^{-1}\right\|\le 1+\alpha.\end{align*}
It follows from \eqref{equ:relationships repeated for n} that for every $n\in\mathbb{N}$,
\begin{align*}\|A_n^\dag\|=&\left\|R_{A_n,N_n}\cdot (A_n)^\dag_{M_nN_n}\cdot L_{A_n,M_n}\right\|\\
\le& \left\|R_{A_n,N_n}\right\|\cdot \left\|(A_n)^\dag_{M_nN_n}\right\|\cdot \left\|L_{A_n,M_n}\right\|\le \alpha (1+\alpha)^2.\end{align*}
This shows that $\sup\limits_n\{\|A_n^\dag\|\}\le \alpha (1+\alpha)^2<+\infty$.
\end{proof}

Our second result on the continuity of the weighted M-P inverse reads as follows.
\begin{theorem}Let $M\in\mathcal{L}(K)$ and $N\in\mathcal{L}(H)$ be weights, and $A\in\mathcal{L}(H,K)$ be given such that $A^\dag_{MN}$ exists. If the sequences $\{M_n\}\subseteq \mathcal{L}(K)_{\mbox{sa}}$ and $\{N_n\}\subseteq \mathcal{L}(H)_{\mbox{sa}}$ satisfy
\begin{equation}\label{approximate with n}\lim_{n\to\infty}\|M_n-M\|=0,\quad \lim_{n\to\infty}\|N_n-N\|=0,\end{equation}
then there exists $n_0\in\mathbb{N}$ such that $A^\dag_{M_nN_n}$ exists for all $n\ge n_0$. Moreover,
$$\lim_{n\to\infty}\|A^\dag_{M_nN_n}-A^\dag_{MN}\|=0.$$
\end{theorem}
\begin{proof} By Theorem~\ref{thm:new formula for the wmp} $A^\dag$ exists, and the operators $R_{A,N}$ and $L_{A,M^{-1}}$ defined by  \eqref{equ:def of R A X} and
\eqref{equ:def of L A Y} are  both invertible. Due to \eqref{approximate with n}, there exists $n_0\in \mathbb{N}$ such that $M_n$, $N_n$  as  well as $R_{A,N_n}$ and $L_{A,M_n^{-1}}$
are all invertible for every $n\ge n_0$, where  $P_{A^*}=A^\dag A$, $P_A=AA^\dag$ and
\begin{equation}\label{defn of Rn and Ln}R_{A,N_n}=P_{A^*}+(I_H-P_{A^*})N_n,\quad L_{A,M_n^{-1}}=P_{A}+M_n^{-1}(I_K-P_{A}).\end{equation}
So by Theorem~\ref{thm:new formula for the wmp}, $A^\dag_{M_nN_n}$ exists for all $n\ge n_0$. Since $R_{A,N_n}^{-1}\to R_{A,N}^{-1}$ and $L_{A,M_n^{-1}}^{-1}\to L_{A,M^{-1}}^{-1}$ as $n\to\infty$, we have
$$\lim_{n\to\infty}A^\dag_{M_nN_n}=\lim_{n\to\infty}R_{A,N_n}^{-1}A^\dag L_{A,M_n^{-1}}^{-1}=R_{A,N}^{-1}A^\dag L_{A,M^{-1}}^{-1}=A^\dag_{MN}.\qedhere$$
\end{proof}

\begin{remark}It is remarkable that there exist certain Hilbert space $H$, and an operator $A\in\mathbb{B}(H)$ together with weights $M,N,M_n,N_n \in\mathbb{B}(H)$ for $n\in \mathbb{N}$
such that
\begin{enumerate}
\item[{\rm (i)}] $A^\dag_{M_nN_n}$ exists for every $n\in\mathbb{N}$;
\item[{\rm (ii)}]$\lim\limits_{n\to\infty}\|M_n-M\|=0$ and $\lim\limits_{n\to\infty}\|N_n-N\|=0$;
\item[{\rm (iii)}]$A^\dag_{MN}$ does not exist.
\end{enumerate}

 Such an example can be provided as follows.
\end{remark}

\begin{example}Given a separable Hilbert space  $E$ and its orthonormal basis $\{e_n:n\in\mathbb{N}\}$, let
  $U\in\mathbb{B}(E)$ be the unilateral shift characterized by
$$Ue_n=e_{n+1}\quad (n\in\mathbb{N}),$$
and $P_1$ be the rank-one projection defined by  $P_1=I_E-UU^*$.  Let
$$H=E\oplus E, \quad M=M_n=I_H,$$ and let  $A,N,N_n\in\mathbb{B}(H)$ be defined by
 \begin{align*} A=\begin{pmatrix}0&U^*\\U&0\end{pmatrix},\quad  N=\begin{pmatrix}0&I_E\\I_E&0\end{pmatrix},\quad N_n=\left(
                      \begin{array}{cc}
                        0 & I_E \\
                        I_E & \frac1n P_1 \\
                      \end{array}
                    \right).
 \end{align*}
 Since
$$\begin{pmatrix}0&I_E\\I_E& X\end{pmatrix}^{-1}=\begin{pmatrix}-X&I_E\\I_E& 0\end{pmatrix},\quad \forall\,X\in\mathbb{B}(E),$$
 we see that $N$ and $N_n(n\in\mathbb{N})$ are weights on $H$.  It is clear that
 $$A^\dag=A,\quad A^\dag A=AA^\dag=\begin{pmatrix}I_E &0\\0& I_E-P_1\end{pmatrix},$$
 so  the operators $R_{A,N_n}$ and $L_{A,M_n^{-1}}$ defined by \eqref{defn of Rn and Ln}  can be calculated as
$$R_{A,N_n}=\left(
        \begin{array}{cc}
          I_E & 0 \\
          P_1 & (I_E-P_1)+\frac{1}{n}P_1\\
        \end{array}
      \right),\quad L_{A,M_n^{-1}}\equiv I_H,$$
which are both invertible in $\mathbb{B}(H)$. Hence, by Theorem~\ref{thm:new formula for the wmp} $A^\dag_{M_nN_n}$ exists for all $n\in\mathbb{N}$.
However, $A^\dag_{MN}$ fails to be existent, since
$$R_{A,N}=\begin{pmatrix}
          I_E & 0 \\
          P_1 & I_E-P_1\end{pmatrix},$$
which is clearly not invertible in $\mathbb{B}(H)$.
 \end{example}

\vspace{2ex}

\noindent\textbf{Declaration of competing interest}

\vspace{2ex}

No known competing financial interests or personal
relationships that could have appeared to influence the work reported in this paper.

\vspace{2ex}

\noindent\textbf{Data availability}

\vspace{2ex}

No data was used for the research described in the paper.

\end{document}